\documentclass[sn-mathphys-num]{sn-jnl}% Math and Physical Sciences Numbered Reference Style
\usepackage{graphicx}%
\usepackage{multirow}%
\usepackage{amsmath,amssymb,amsfonts}%
\usepackage{amsthm}%
\usepackage{mathrsfs}%
\usepackage[title]{appendix}%
\usepackage{xcolor}%
\usepackage{textcomp}%
\usepackage{manyfoot}%
\usepackage{booktabs}%
\usepackage{algorithm}%
\usepackage{algorithmicx}%
\usepackage{algpseudocode}%
\usepackage{listings}%
\theoremstyle{thmstyleone}%
\newtheorem{theorem}{Theorem}%  meant for continuous numbers
\newtheorem{proposition}[theorem]{Proposition}%

\theoremstyle{thmstyletwo}%
\newtheorem{remark}{Remark}%

\theoremstyle{thmstylethree}%
\newtheorem{lemma}[theorem]{Lemma}%

\newcommand{\R}{\mathbb{R}}
\def\E{\mathrm{E}}
\def\Var{\mathrm{Var}}

\def\d{\mathrm{d}}
\def\Cov{\mathrm{Cov}}
\def\R{\mathbb{R}}

\def\e{\varepsilon}
\def\b{\beta}
\def\u{\upsilon}
\def\p{\boldsymbol p}
\begin{document}

\title[Functional CLTs for parabolic Anderson model]{Functional central limit theorems for spatial averages of the parabolic Anderson model with delta initial condition in dimension $d\geq 1$}


\author[1]{\fnm{Wanying} \sur{Zhang}}\email{wyzhang20@mails.jlu.edu.cn}

\author*[1]{\fnm{Yong} \sur{Zhang}}\email{zyong2661@jlu.edu.cn}

\author[1]{\fnm{Jingyu} \sur{Li}}\email{lijingyu@jlu.edu.cn}

\affil*[1]{\orgdiv{School of Mathematics}, \orgname{Jilin University}, \orgaddress{\city{Changchun}, \postcode{130012}, \state{Jilin}, \country{China}}}

%\affil[2]{\orgdiv{Department}, \orgname{Organization}, \orgaddress{\street{Street}, \city{City}, \postcode{10587}, \state{State}, \country{Country}}}

%\affil[3]{\orgdiv{Department}, \orgname{Organization}, \orgaddress{\street{Street}, \city{City}, \postcode{610101}, \state{State}, \country{Country}}}


\abstract{Let $\{u(t,x)\}_{t>0,x\in{{\mathbb R}^{d}}}$ denote the solution to a $d$-dimensional parabolic Anderson model with delta initial condition
 and driven by a multiplicative noise that is white in time and has a spatially homogeneous covariance given by a nonnegative-definite measure $f$. Let
$S_{N,t}:=N^{-d}\int_{{[0,N]}^d}{[U(t,x)-1]}{\rm d}x$ denote the spatial average on
${{\mathbb R}^{d}}$. We obtain various functional central limit theorems (CLTs) for spatial averages based on the quantitative analysis of $f$ and spatial dimension $d$. In particular, when $f$ is given by Riesz kernel, that is, $f({\rm x})={\Vert x \Vert}^{-\beta}{\rm d}x$, $\b\in(0,2\wedge d)$, the functional CLT is also based on the index $\beta$.}



%%\pacs[MSC Classification]{60F17, 60H15, 60H07}

\maketitle

\section{Introduction}\label{sec1}

Consider the following parabolic Anderson model:
{\begin{equation}\label{1.1}
\left\{
\begin{array}{lr}
{\partial _t}{u(t,x)}={\frac{1}{2}}{\Delta u(t,x)}+{u(t,x)}{\eta (t,x)}~~~{\rm for} ~{(t,x)} \in {(0,+\infty)}{\times{\mathbb R}^{d}},\\
{\rm subject~to}~~~~~{u(0)=\delta_{0}},\\
\end{array}
\right.
\end{equation}
where $\eta$ denotes a centered, generalized Gaussian random field such that for all $s,t\geq0$ and $x,y \in {\mathbb R}^{d}$
\begin{align}\label{1.2}
{\rm Cov}[\eta(t,x),\eta(s,y)]=\delta_{0}(t-s)f(x-y)~{\rm for~all}~ s,t\geq0 ~{\rm and} ~x,y \in {\mathbb R}^{d},
\end{align}
for a non-zero, nonnegative-definite, tempered Borel measure $f$ on ${\mathbb R}^{d}$. In order to avoid degeneracies, we assume that $f({\mathbb R}^{d})>0$.}

{Following from \cite{MR0876085}, the equation $(\ref{1.1})$ has a mild solution that satisfies the integral equation:
\begin{align}\label{1.3}
u(t,x)={\boldsymbol p}_{t}(x)+\int_{(0,t)\times{\mathbb R}^{d}}{\boldsymbol p}_{t-s}(x-y)u(s,y)\eta({\rm d}s,{\rm d}y)~~~{\rm for}~t>0,~{\rm and}~x\in {\mathbb R}^{d},
\end{align}
where ${\boldsymbol p}_{t}(x)$ is the heat kernel defined by
$$
{\boldsymbol p}_{t}(x)={{(2\pi t)}}^{-d/2}e^{-\Vert x\Vert^2/2t}~~~{\rm for}~t>0,~{\rm and}~x\in{\mathbb R}^{d}.
$$
Because
$$
\frac{{\boldsymbol p}_{t-s}(a){\boldsymbol p}_{s}(b)}{{\boldsymbol p}_{t}(a+b)}={\boldsymbol p}_{s(t-s)/t}\left(b-\frac{s}{t}(a+b)\right)~~~{\rm for~all}~0<s<t,~{\rm and}~a,b\in{\mathbb R}^{d},
$$
(\ref{1.3}) can be recast as the following evolution equation:
\begin{equation}\label{1.4}
U(t,x)=1+\int_{(0,t)\times{\mathbb R}^{d}}{\boldsymbol p}_{s(t-s)/t}\left(y-\frac{s}{t}x\right)U(s,y)\eta({\rm d}s,{\rm d}y),
\end{equation}
where $U(t,x):=u(t,x)/{\boldsymbol p}_t(x).$}

{In order to ensure the existence and uniqueness of $u$, hence also $U$, we assume that the Fourier transform ${\hat f}$ satisfies the integrability condition:
$$
\Upsilon(\beta):=\frac{1}{{(2\pi)}^{d}}\int_{{\mathbb R}^{d}}\frac{{\hat f({\rm d}y)}}{\beta+{\Vert y \Vert}^{2}}<\infty~~~{\rm for~all}~\beta>0.
$$}

The purpose of this paper is to establish specific functional CLTs for spatial averages of $U$. The reason we do not consider the spatial average of $u$ is that $\{u(t,x)\}_{x\in\R}$ lacks spatial stationarity, whereas the renormalized process $\{U(t,x)\}_{x\in \R}$ is spatial stationary \cite{MR4334682,MR4242879}. Let us start by introducing the spatial average:
\begin{equation}\label{1.5}
S_{N,t}=\frac{1}{N^{d}}\int_{{[0,N]}^d}{[U(t,x)-1]}{\rm d}x.
\end{equation}
In order to present the main theorems, we introduce a quantity associated with $f$. For every real number $m>0$, define
\begin{equation}\label{1.6}
I_{m}(x):=m^{-d}\mathbf{1}_{{[0,m]}^{d}}(x),~~~{\widetilde{I}}_{m}(x):=I_{m}(-x)~~~{\rm for}~x\in{\mathbb R}^{d}
\end{equation}
and
\begin{align}\label{1.7}
{\mathcal{R}}
(f):&=\frac{1}{{\pi}^{d}}\int_{0}^{\infty}{\rm d}s\int_{{\mathbb R}^{d}}{\hat f({\rm d}z)}\prod_{j=1}^{d}\frac{1-\cos(sz_{j})}{{(sz_{j})}^{2}}\notag\\
&=\frac{1}{{(2\pi)}^{d}}\int_{0}^{\infty}\frac{{\rm d}r}{{r}^d}\int_{{\mathbb R}^{d}}{\hat f({\rm d}z)}\left(\widehat{{\left(I_{1}\ast {\widetilde{I}}_{1}\right)(\bullet/r)}}\right)(z).
\end{align}
\begin{theorem}\label{thm1.1}
 Assume ${\mathcal R}(f)<\infty$, then for any fixed real number $T>0$, we have
\begin{equation}\label{1.8}
\{\sqrt{N}S_{N,t}\}_{t\in[0,T]}{\xrightarrow{C[0,T]}}\{\mathcal{G}_{t}\}_{t\in[0,T]}~~~{\rm{as}}~N\rightarrow\infty,
\end{equation}
where $``{\xrightarrow{C[0,T]}}"$ denotes weak convergence on the space of all continuous functions $C[0,T]$, $\{\mathcal{G}_{t}\}_{t\geq0}$ is a centered Gaussian process with covariance $\E[\mathcal{G}_{t_{1}}\mathcal{G}_{t_{2}}]={ g}_{t_{1},t_{2}}$, which is defined explicitly in Proposition $\ref{prop3.1}$.
\end{theorem}

\begin{theorem}\label{thm1.2}
 Assume $d=1$ and $f$ is a Rajchman measure, that is, $f$ satisfies $f({\mathbb{R}})<\infty$ and $\lim_{{x}\rightarrow\infty}\hat{f}(x)=0$,
then for any fixed real number $T>0$, we have
\begin{equation}\label{1.10}
\left\{\sqrt{\frac{N}{\log N}}S_{N,t}\right\}_{t\in[0,T]}{\xrightarrow{C[0,T]}}\left\{\sqrt{f(\mathbb{R})}B_{t}\right\}_{t\in[0,T]},
\end{equation}
where $B$ denotes a standard Brownian motion.
\end{theorem}
\begin{remark}
The collection of Rajchman measures \cite{MR1364897} is very rich. For example, $f$ is given by a Gaussian kernel ($f(\d x)={\boldsymbol p}_1(x)\d x$ and ${\hat f}(\d x)=e^{-\vert x \vert^2/2}\d x$) or a Cauchy kernel ($f(\d x)=(1+\vert x \vert^2)^{-1}\d x$ and ${\hat f}(\d x)=e^{-\vert x \vert}\d x$).
\end{remark}
\begin{remark}
In fact, the assumption $\lim_{{x}\rightarrow\infty}\hat{f}(x)=0$ is not necessary. We can obtain the functional CLT provided the limit of the covariance exists. It is worth noting that the limit can only take values between $(t_{1}\wedge t_{2})f(\mathbb R)$ and $2(t_{1}\wedge t_{2})f(\mathbb R)$; see Proposition \ref{prop3.2}. \vskip 0.8em
\end{remark}
\begin{remark}
 One can see from Lemma 5.9 in \cite{MR4242879} that $\mathcal{R}(f)=\infty$ if $d=1$. Therefore, Theorem \ref{thm1.1} holds only for $d\geq2$. And there does not seem to be a canonical functional CLT when $\mathcal{R}(f)=\infty~(d\geq2)$, or $f(\mathbb R)=\infty$.  In this situation, we present a special case that $f$ is given by Riesz kernel as an example.
\end{remark}
\begin{theorem}\label{thm1.3}$\rm (Riesz~kernel).$
 Assume $f({\rm d}x)={\Vert x\Vert}^{-\beta}{\rm d}x,$ for some $\beta\in(0,d\wedge2)$, then for any fixed real number $T>0$, we have\\
~\\
${{\mathbf {(A)}}}$ If $\beta\in(0,1)$, then\\
\begin{equation}\label{1.11} \left\{\sqrt{N^{\beta}}S_{N,t}\right\}_{t\in[0,T]}{\xrightarrow{C[0,T]}}\left\{{\mathcal C}_{t}^{(1)}\right\}_{t\in[0,T]}~~~{\rm as}~N\rightarrow\infty.
\end{equation}
~\\
${{\mathbf {(B)}}}$ If $\beta=1$, then\\
\begin{equation}\label{1.12}
 \left\{\sqrt{\frac{N}{\log N}}S_{N,t}\right\}_{t\in[0,T]}{\xrightarrow{C[0,T]}}\left\{{\mathcal C}_{t}^{(2)}\right\}_{t\in[0,T]}~~~{\rm as}~N\rightarrow\infty.
\end{equation}
~\\
${{\mathbf {(C)}}}$ If $\beta\in(1,2\wedge d)$, then\\
\begin{equation}\label{1.13}
 \left\{\sqrt{N^{2-\beta}}S_{N,t}\right\}_{t\in[0,T]}{\xrightarrow{C[0,T]}}\left\{{\mathcal C}_{t}^{(3)}\right\}_{t\in[0,T]}~~~{\rm as}~N\rightarrow\infty,
\end{equation}
where \small{$\left\{{\mathcal C}_{t}^{(i)}\right\}_{t\geq0}$}$(i=1,2,3)$ is a centered Gaussian process with \small{$\E[{{{\mathcal C}_{t_{1}}^{(i)}{{\mathcal C}_{t_{2}}^{(i)}}}}]={ c}^{(i)}_{t_{1},t_{2}}$}, which is defined explicitly in Proposition $\ref{prop3.3}$.
\end{theorem}
There are many arguments for CLTs and functional CLTs for spatial averages of the solution to (\ref{1.1}). In the case of constant/flat initial condition $u(0)\equiv1$, the quantitative CLT and the related functional CLT were introduced first in \cite{MR4167203} by using Malliavin-Stein's method for $d=1,~f=\delta_{0}$. As for $d\geq1$, the condition $f({\mathbb R}^{d})<\infty$ implies a standard form of CLT \cite[Theorem 1.1]{MR4421618}, and \cite{MR4563698} deduced the convergence rate for the CLT in terms of total variation distance and established a corresponding functional CLT. As an example of the condition $f({\mathbb R}^{d})=\infty$, the relative discussions in the case that $f={\Vert x \Vert}^{-\beta}$ ($0<\beta<2\wedge d$) were given by \cite{MR4098872}. While in the case of delta initial condition $u(0)=\delta_{0}$, \cite{MR4334682} proved that the convergence rate for the variance of the spatial average behaves differently from the rate in the case $u(0)\equiv1$ when $d=1$ and $f=\delta(0)$, this also affects the form of the functional CLT. Later, \cite{MR4242879} generalized the CLT to the multidimensional case, which is much more involved than the case $d=1$. The results depend on not only the dimension $d$ but also the behavior of $f$. CLT and its related variations for other types of equations can be found in \cite{MR4337703,MR4441503,MR4683383,MR4479916}. For other limit theories, we refer to \cite{MR3949959,MR4066510,MR4535893,MR4578556,DOI:10.1080/03610926.2024.2335541,MR3620268}.
%Other relevant CLT and its variations can be observed in many literatures \cite{Jing21,Nualart21,Nualart221}.

In this paper, we establish the functional CLTs in detail under the condition $d\geq1$  mainly based on Malliavin-Stein's method, Fourier analysis and Poincar$\acute{\rm e}$-type inequality. The Functional CLT is an extension of the CLT to function spaces, providing a more comprehensive understanding of the limiting behavior of $\{S_{N,t}\}_{t>0}$ as a stochastic process. Although the methods for proving functional central limit theorems is already well-developed \cite{MR4479916,MR4167203,MR4098872,MR4334682,MR4563698,MR4683383,MR4242879}, it is undeniable that the $\delta_0$ initial condition and the measure $f$ introduce significant challenges to our proof. We derive the asymptotic covariance of the spatial averages in Section 3, which is an extension on both the asymptotic behavior of the variance \cite[Section 5]{MR4242879} and the asymptotic behavior of the covariance for $f=\delta_0$ \cite[Section 4]{MR4334682}. In Section 4, we prove the uniform upper bound and tightness of the spatial averages. In order to prove our main theorems, it remains to establish the convergence of finite dimensional distributions, which is shown in Section 5. And the last section is Appendix, we introduce a few technical lemmas that are used throughout the paper.

Throughout this paper, we write ``$g_{1}(x)\lesssim g_{2}(x)$ for all $x\in X$'' when there exists a real number $L$ independent of $x$ such that $g_{1}(x)\leq Lg_{2}(x)$ for all $x\in X$. And for every $Z\in L^{k}(\Omega)$, we write ${\Vert Z \Vert}_{k}$ instead of ${\Vert Z \Vert}_{L^{k}(\Omega)}$.

\section{Preliminaries}\label{sec2}
Let us begin by introducing the fundamentals of Malliavin calculus and stochastic integral. Let ${\mathcal H}_{0}$ denote the reproducing kernel Hilbert space, spanned by all real-valued functions on ${\mathbb R}^{d}$, that corresponds to the scalar product $\langle\phi,\psi\rangle:={\langle\phi,{\psi\ast f\rangle}}_{L^{2}({\mathbb R}^{d})}$, and let ${\mathcal H}=L^{2}({\mathbb R}_{+}\times{\mathcal H}_{0})$. The Gaussian random field $\{W(h)\}_{h\in{\mathcal H}}$ formed by such Wiener integrals
\begin{equation}\label{2.1}
W(h):=\int_{{\mathbb R}^{+}\times{\mathbb R}^{d}}h(s,y)\eta({{\rm d}s,{\rm d}y})
\end{equation}
defines an isonormal Gaussian process on the Hilbert space ${\mathcal H}$. Therefore, we can develop the Malliavin calculus (see, for instance, \cite{MR2200233}). We denote by $D$ the Malliavin derivative operator and by $\delta$ the corresponding divergence operator, then we have the following important property:

Let $F$ denote a predictable and square-integrable random field valued in the Gaussian Sobolev space ${\mathbb D}^{1,2}$, and also, $F$ belongs to the domain of $\delta$ denoted by Dom[$\delta$], then $\delta(F)$ coincides with the Walsh integral:
$$
\delta(F)=\int_{{\mathbb R}^{+}\times{\mathbb R}^{d}}F(s,y)\eta({{\rm d}s,{\rm d}y}).
$$
Moreover, for all $p\geq2$, there exists $c_p>0$, such that for all $t>0$, the following Burkholder-Davis-Gundy (BDG) inequality holds:
\begin{align}
\E\left[\left\vert\int_{(0,t)\times\R^d}F(s,y)\eta(ds,dy)\right\vert^p\right]\leq c_p\E\left[\left(\int_{(0,t)\times\R^{2d}}F(s,y)F(s,y+y^\prime)f(\d y')\d y\d s\right)^{\frac{p}{2}}\right]\notag
\end{align}

Let us recall the following inequalities from Theorem 1.1 and Proposition 4.1 in \cite{MR4242879}.

\begin{lemma}\label{lem2.1}
Suppose $U={\{U(t,x)\}}_{t>0,x\in{\mathbb R}^{d}}$ is a predictable solution to the integral equation $(\ref{1.4})$, then, $U$ is the only
predictable solution to $(\ref{1.4})$ that satisfies the following for all ${\varepsilon}\in(0,1)$, $t>0$, and $k\geq2$:
\begin{equation}\label{2.2}
\sup_{x\in{\mathbb R}^{d}}{\Vert U(t,x) \Vert}_{k}\leq\left(\frac{2}{\varepsilon}\right)\exp\left\{{\frac{4}{t}\Upsilon^{-1}\left(\frac{1-\varepsilon}{4z^{2}_{k}}\right)}\right\}:=c_{t,k},
\end{equation}
where $z_{k}$ denotes the optimal constant in the BDG inequality. And $\{U(t,x)\}_{x\in{\mathbb R}^d}$ is a  stationary random field for every $t>0$. Moreover, for almost every $(s,y)\in(0,t)\times{\mathbb R}^{d}$, there exists a real number $C_{t,k}$ such that
\begin{equation}\label{2.3}
{\Vert D_{s,y}U(t,x)\Vert}_{k}\leq C_{t,k}{\boldsymbol p}_{s(t-s)/t}\left(y-\frac{s}{t}x\right).
\end{equation}
\end{lemma}
The following lemma, which is a generalization of a result of Theorem 6.1.2 in \cite{MR2962301}, plays an important role in the proofs of our theorems.
\begin{lemma}\label{lem2.2}
Let $F=(F^{(1)},\ldots,F^{(m)})$ denote a random vector such that for every $i=1,\ldots,m$, $F^{(i)}=\delta(v^{(i)})$ for some $v^{(i)}\in {\rm Dom}[\delta]$. Assume additionally that $F^{(i)}\in{\mathbb D}^{1,2}$ for $i=1,\ldots,m$. Let $G$ be a centered m-dimensional Gaussian random vector with covariance matrix $(C_{i,j})_{1\leq i,j\leq m}$. Then, for any $h\in C^{2}({\mathbb R}^{m})$ that has bounded partial derivatives, we have
$$
\left\vert \E (h(F))-\E (h(G))\right\vert \leq \frac{1}{2}{\Vert h^{''}\Vert}_{\infty}\sqrt{\sum_{i,j=1}^{m}\E{\left({\left\vert C_{i,j}-{\langle DF^{(i)},v^{(j)}\rangle}_{\mathcal H}\right\vert}^{2}\right)}},
$$
where
$$
{\Vert h^{''}\Vert}_{\infty}:=\max_{1\leq i,j\leq m}\sup_{x\in{\mathbb R}^{m}}\left\vert\frac{\partial^{2}h(x)}{\partial x_{i} \partial x_{j}}\right\vert.
$$
\end{lemma}

Finally, we recall the following  Poincar\rm{$\acute{\rm e}$}-type inequality \cite[Section 2.2]{MR4242879}:
\begin{lemma}Suppose that $F,G\in {\mathbb D}^{1,2}$ and $DF$ and $DG$ are real-valued random variables, then we have
\begin{equation}\label{2.4}
\vert{\rm Cov}(F,G)\vert\leq \int_{0}^{\infty}{\rm d}s\int_{{\mathbb R}^{d}}{\rm d}y\int_{{\mathbb R}^{d}}f({\rm d}y'){\Vert D_{s,y}F\Vert}_{2}{\Vert D_{s,y+y'}G\Vert}_{2},
\end{equation}
\end{lemma}

\section{Asymptotic behavior of the covariance}
Recall that the spatial average $S_{N,t}$ and the quantity ${\mathcal R}(f)$ are defined respectively in (\ref{1.5}) and (\ref{1.7}). In this section, we will estimate the asymptotic behavior of the covariance functions of $S_{N,t}$.

In order to simplify the exposition, for any $t_{1},t_{2}>0$ we define
\begin{equation}\label{3.1}
\tau:=\frac{2t_{1}t_{2}}{t_{1}+t_{2}},~\tau_{1}:=\frac{2t_{2}}{t_{1}+t_{2}},~{\rm and}~~\tau_{2}:=\frac{2t_{1}}{t_{1}+t_{2}}.
\end{equation}

\begin{proposition}\label{prop3.1}
$(d\geq2).$ Assume ${\mathcal R}(f)<\infty$, then for any $t_{1}, t_{2}>0,$ we have
\begin{equation}\label{3.2}
\lim_{N\rightarrow\infty}\Cov(\sqrt{N}S_{N,t_{1}},\sqrt{N}S_{N,t_{2}})=g_{t_{1},t_{2}},
\end{equation}
where
\begin{equation*}
g_{t_{1},t_{2}}=\frac{\tau}{{(2\pi)}^{d}}\int_{0}^{\infty}\frac{\rm ds}{s^d}\int_{{\mathbb R}^{d}}\left[\widehat{\left(I_{({\tau}_{1}\wedge {\tau}_{2})}\ast \tilde{I}_{(\tau_{1}\vee\tau_{2})}\right)\left(\frac{\bullet}{s}\right)}\right](z){\hat f}({\rm d}z).
\end{equation*}
\end{proposition}

\begin{proposition}\label{prop3.2}
$(d=1)$. Assume $f(\mathbb R)<\infty$, then for any $t_{1},t_{2}>0$, we have
\begin{equation}\label{3.3}
\begin{aligned}
&\liminf_{N\rightarrow\infty}\Cov\left(\sqrt{\frac{N}{\log N}}S_{N,t_{1}},\sqrt{\frac{N}{\log N}}S_{N,t_{2}}\right)\geq(t_{1}\wedge t_{2})f(\mathbb R),\\
&\limsup_{N\rightarrow\infty}\Cov\left(\sqrt{\frac{N}{\log N}}S_{N,t_{1}},\sqrt{\frac{N}{\log N}}S_{N,t_{2}}\right)\leq 2(t_{1}\wedge t_{2})f(\mathbb R).
\end{aligned}
\end{equation}
In particular, if $\lim_{x\rightarrow\infty}{\hat f}(x)=0$, then the limit of the covariance exists, we have that
\begin{equation}\label{3.4}
\begin{aligned}
\lim_{N\rightarrow\infty}\Cov\left(\sqrt{\frac{N}{\log N}}S_{N,t_{1}},\sqrt{\frac{N}{\log N}}S_{N,t_{2}}\right)= (t_{1}\wedge t_{2})f(\mathbb R).
\end{aligned}
\end{equation}
\end{proposition}

\begin{proposition}\label{prop3.3}
{\rm {(Riesz kernel)}}.
Assume $f(dx)={\Vert x\Vert}^{-\beta}{\rm d}x$ and ${\hat f}(dx)={\kappa}_{\beta,d}{\Vert x\Vert}^{\beta-d}{\rm d}x$ for $0<\beta<2\wedge d$, where ${\kappa}_{\beta,d}$ is a real number depending on $\beta$ and $d$, we have

${{\mathbf {(A)}}}$ If $\beta\in(0,1)$, then\\
\begin{equation}\label{3.5}
\lim_{N\rightarrow \infty}\Cov\left(\sqrt{N^{\beta}}S_{N,t_{1}},\sqrt{N^{\beta}}S_{N,t_{2}}\right)=c^{(1)}_{t_{1},t_{2}},
\end{equation}
where $$
c^{(1)}_{t_{1},t_{2}}=\frac{{(t_{1}\wedge t_{2})^{1-\beta}}{\tau^{\beta}}}{1-\beta}\int_{{\mathbb R}^d}{\left(I_{({\tau}_{1}\wedge {\tau}_{2})}\ast \tilde{I}_{(\tau_{1}\vee\tau_{2})}\right)}(z){\Vert z\Vert}^{-\beta}{\rm d}z.
$$

${{\mathbf {(B)}}}$ If $\beta=1$, then\\
\begin{equation}\label{3.6}
\lim_{N\rightarrow \infty}\Cov\left(\sqrt{\frac{N}{\log N}}S_{N,t_{1}},\sqrt{\frac{N}{\log N}}S_{N,t_{2}}\right)=c^{(2)}_{t_{1},t_{2}},
\end{equation}
where $$
c^{(2)}_{t_{1},t_{2}}=\frac{2\tau \kappa_{1,d}}{{(2\pi)}^{d}}\int_{{\mathbb R}^d}\widehat{{\left(I_{({\tau}_{1}\wedge {\tau}_{2})}\ast \tilde{I}_{(\tau_{1}\vee\tau_{2})}\right)}}(z){\Vert z \Vert}^{1-d}{\rm d}z.
$$

${{\mathbf {(C)}}}$ If $\beta\in(1,2\wedge d)$, then\\
\begin{equation}\label{3.7}
\lim_{N\rightarrow \infty}\Cov\left(\sqrt{N^{2-\beta}}S_{N,t_{1}},\sqrt{N^{2-\beta}}S_{N,t_{2}}\right)=c^{(3)}_{t_{1},t_{2}},
\end{equation}
where $$
c^{(3)}_{t_{1},t_{2}}=\frac{{\tau^{2-\beta}\kappa_{\beta,d}}}{{(2\pi)}^d}\int_{{\mathbb R}^d}\widehat{{\left(I_{({\tau}_{1}\wedge {\tau}_{2})}\ast \tilde{I}_{(\tau_{1}\vee\tau_{2})}\right)}}(z){\Vert z \Vert}^{2-\beta-d}{\rm d}z\int_{0}^\infty{r^{\beta-2}}e^{-r}{\rm d}r.
$$
\end{proposition}

\begin{remark}
\rm {Unlike the case that $u(0)\equiv1$, we can't obtain the asymptotic behavior of the covariance from the asymptotic variance directly due to the complexity of the condition probability density $\p_{s(t-s)/t}(sx/t)$. Therefore, we need to prove the above propositions as follows. }
\end{remark}
For every $y\in{\mathbb R}^{d}$, and $s>0$, we first denote
\begin{equation}\label{3.8}
\chi_{s}(y):=\Cov[U(s,0),U(s,y)]=\E[U(s,0)U(s,y)]-1.
\end{equation}
By (\ref{1.4}), the stationary property of $\{U(t,x)\}_{x\in{\mathbb R}^{d}}$ and the semigroup property of the heat kernel,
\begin{align*}
\Cov(S_{N,{t_{1}}},S_{N,{t_{2}}})
&=\frac{1}{N^{2d}}\int_{[0,N]^{2d}}{\rm d}x_{1}{\rm d}x_{2}~\Cov[U(t_{1},x_{1})-1, U(t_{2},x_{2})-1]\\
&=\frac{1}{N^{2d}}\int_{[0,N]^{2d}}{\rm d}x_{1}{\rm d}x_{2}\int_{(0,t_{1}\wedge t_{2})\times{{\mathbb R}^{2d}}}f({\rm d}y'){\rm d}y{\rm d}s~\E[U(s,y)U(s,y+y')]
\end{align*}
\begin{align*}
&\times {\boldsymbol p}_{{s[(t_{1}\wedge t_{2})-s]}/(t_{1}\wedge t_{2})}\left( y-\frac{s}{t_{1}\wedge t_{2}}x_{1}\right){\boldsymbol p}_{{s[(t_{1}\vee t_{2})-s]}/(t_{1}\vee t_{2})}\left( y+y'-\frac{s}{t_{1}\vee t_{2}}x_{2}\right)\\
=&\frac{1}{N^{2d}}\int_{[0,N]^{2d}}{\rm d}x_{1}{\rm d}x_{2}\int_{(0,t_{1}\wedge t_{2})\times{{\mathbb R}^{d}}}f({\rm d}y'){\rm d}s~(1+\chi_{s}(y'))\\
&\times {\boldsymbol p}_{{s(t_{1}-s)}/t_{1}+{s(t_{2}-s)}/t_{2}} \left(y'-\left(\frac{s}{t_{1}\vee t_{2}}x_{2}-\frac{s}{t_{1}\wedge t_{2}}x_{1}\right)\right).
\end{align*}
Suppose now that $t_{1}<t_{2}$, then we can write
\begin{align*}
\Cov(S_{N,{t_{1}}},S_{N,{t_{2}}})
&=V^{(1)}_{N,t_1,t_2}+V^{(2)}_{N,t_1,t_2},
\end{align*}
where
\begin{align}
V^{(1)}_{N,t_1,t_2}
&=\frac{1}{N^{2d}}\int_{[0,N]^{2d}}{\rm d}x_{1}{\rm d}x_{2}\int_{0}^{t_{1}}{\rm d}s\int_{{\mathbb R}^{d}}f({\rm d}y')~{\boldsymbol p}_{{2s(\tau-s)}/\tau}\left(y'-\frac{s}{\tau}
(\tau_{2}x_{2}-\tau_{1}x_{1})\right),
\label{3.9}\\
V^{(2)}_{N,t_1,t_2}
&=\frac{1}{N^{2d}}\int_{[0,N]^{2d}}{\rm d}x_{1}{\rm d}x_{2}\int_{0}^{t_{1}}{\rm d}s\int_{{\mathbb R}^{d}}f({\rm d}y')~\notag\\
&~~~~~~~~~~~~~~~~~~~~~~~~~~~~~~~~~~\times{\boldsymbol p}_{{2s(\tau-s)}/\tau}\left(y'-\frac{s}{\tau}
(\tau_{2}x_{2}-\tau_{1}x_{1})\right)\chi_{s}(y').
\label{3.10}
\end{align}
Furthermore, by making the change of variables $[x_i\mapsto Nx_i/\tau_i,~{\rm for}~i=1,2]$, since $\chi_{t}(x)\leq \chi_{t}(0)$ for all $t>0$ and $x\in{{\mathbb R}^{d}}$ \cite[Lemma 5.6]{MR4242879}, we have
\begin{align}
V^{(1)}_{N,t_1,t_2}
&=\int_{{{\mathbb R}^{d}}}{\rm d}z\left(I_{\tau_2}\ast {\tilde I_{\tau_1}}\right)(z)
\int_{0}^{t_{1}}{\rm d}s~\left({\boldsymbol p}_{{2s(\tau-s)}/\tau}\ast f\right)\left(\frac{sN}{\tau}z\right),\label{3.11}\\
V^{(2)}_{N,t_1,t_2}
&\leq\int_{{{\mathbb R}^{d}}}{\rm d}z\left(I_{\tau_2}\ast {\tilde I_{\tau_1}}\right)(z)\int_{0}^{t_1}{\rm d}s~\chi_{s}(0)\left({\boldsymbol p}_{{2s(\tau-s)}/\tau}\ast f\right)\left(\frac{sN}{\tau}z\right),\label{3.12}
\end{align}
where $\tau,\tau_{1},\tau_{2}$ are defined in (\ref{3.1}). Next, we define
\begin{equation}\label{3.13}
\psi_{\tau_1,\tau_2}(z):=\widehat{\left(I_{\tau_2}\ast {\tilde I_{\tau_1}}\right)}(z)=\widehat{{I}_{\tau_2}}(z)\overline{\widehat{{{{ I}}_{\tau_1}}}(z)}~~~{\rm for~all~}z\in\R^d.
\end{equation}
From Lemma \ref{lem6.1} (1), we have
\begin{equation}\label{3.13'}
\vert \psi_{\tau_1,\tau_2}(z)\vert\lesssim1\wedge\Vert z\Vert^{-2},
\end{equation}
which we use several times in computing integrals.
Following the argument for $V^{(2)}_{N,t_1,t_2}$ in Section 5 of \cite{MR4242879},
we can see that the main asymptotic behavior of the covariance is $V^{(1)}_{N,t_1,t_2}$ by replacing $t_1$ with $\tau$ in (\ref{3.12}) and using the inequality (\ref{3.13'}). Therefore, it remains to consider the behavior of $V^{(1)}_{N,t_1,t_2}$ only.

\subsection{Analysis in $d\geq2$}
\begin{proof}[Proof of Proposition \ref{prop3.1}] By Fubini's theorem and a change of variables, we can see that
\begin{equation*}
\begin{aligned}
V^{(1)}_{N,t_1,t_2}
&=\frac{\tau}{N}\int_{0}^{{t_{1}N}/{\tau}}{\rm d}s\int_{{{\mathbb R}^{d}}}\left(I_{\tau_2}\ast {\tilde I_{\tau_1}}\right)(z)\left({\boldsymbol p}_{{\frac{2s\tau}{N}(1-\frac{s}{N})}}\ast f\right)\left(sz\right) \\
&=\frac{\tau}{N}\int_{0}^{{t_{1}N}/{\tau}}\frac{{\rm d}s}{s^d}\int_{{{\mathbb R}^{d}}}\left(I_{\tau_2}\ast {\tilde I_{\tau_1}}\right)\left(\frac{z}{s}\right)\left({\boldsymbol p}_{{\frac{2s\tau}{N}(1-\frac{s}{N})}}\ast f\right)\left(z\right) \\
&=\frac{\tau}{{(2\pi)}^{d}N}\int_{0}^{{t_{1}N}/{\tau}}\frac{{\rm d}s}{s^d}\int_{{{\mathbb R}^{d}}}\hat f({\rm d}z)\widehat{\left[\left(I_{\tau_2}\ast {\tilde I_{\tau_1}}\right)\left(\frac{\bullet}{s}\right)\right]}(z)
\exp\left[{{-\frac{s\tau}{N}\left(1-\frac{s}{N}\right)}}\right].
\end{aligned}
\end{equation*}
 Then we apply Lemma \ref{lem6.1} (2) and the dominated convergence theorem to find that
\begin{equation*}
\begin{aligned}
\lim_{N\rightarrow\infty}NV^{(1)}_{N,t_1,t_2}
=\frac{\tau}{{(2\pi)}^{d}}\int_{0}^{\infty}\frac{{\rm d}s}{s^d}\int_{{{\mathbb R}^{d}}}\hat f({\rm d}z)\widehat{\left[\left(I_{\tau_2}\ast {\tilde I_{\tau_1}}\right)\left(\frac{\bullet}{s}\right)\right]}(z).
\end{aligned}
\end{equation*}
This implies (\ref{3.2}).
\end{proof}

\subsection{Analysis in $d=1$}
\begin{lemma}\label{lem3.6}
For all $t_{1},t_{2}>0$,
\begin{equation}\label{3.15''}
\begin{aligned}
(t_{1}\wedge t_{2})f(\mathbb R)\leq\liminf_{N\rightarrow\infty}\frac{N}{\log N}V^{(1)}_{N,t_1,t_2}\leq\limsup_{N\rightarrow\infty}\frac{N}{\log N}V^{(1)}_{N,t_1,t_2}\leq 2(t_{1}\wedge t_{2})f(\mathbb R).
\end{aligned}
\end{equation}
\end{lemma}
\begin{proof}
Recalling (\ref{3.11}), we can write
\begin{align*}
V^{(1)}_{N,t_1,t_2}
&=\frac{\tau }{2\pi N}\int_{0}^{t_{1}}\frac{{\rm d}s}{s}\int_{-\infty}^{\infty}{\rm d}z~\psi_{\tau_1,\tau_2}(z)\exp\left(-\frac{\tau(\tau-s)}{N^2s}z^2\right){\hat f}\left(\frac{\tau z}{Ns}\right)\\
&\leq\frac{\tau}{2\pi N}\int_{0}^{\tau}\frac{{\rm d}s}{s}\int_{-\infty}^{\infty}{\rm d}z~\psi_{\tau_1,\tau_2}(z)\exp\left(-\frac{\tau(\tau-s)}{N^2s}z^2\right){\hat f}\left(\frac{\tau z}{Ns}\right).
\end{align*}
Then, by using Lemma \ref{lem6.4}, the dominated convergence theorem and $\frac{1}{2\pi}\int_{-\infty}^\infty \psi_{\tau_1,\tau_2}(z)\d z=1/\tau_1$, we obtain the third inequality in (\ref{3.15''}).

Next, we prove the first inequality. With a change of variable $s=t_{1}rN^{-2}$, the term $V^{(1)}_{N,t_1,t_2}$ can be recast as
\begin{align*}
V^{(1)}_{N,t_1,t_2}
&=\frac{\tau }{2\pi N}\int_{0}^{N^{2}}\frac{{\rm d}r}{r}\int_{-\infty}^{\infty}{\rm d}z~\psi_{\tau_1,\tau_2}(z)\exp\left[-\tau z^2\left(\frac{1-\frac{r}{\tau_{1}N^2}}{{r}/{\tau_1}}\right)\right]{\hat f}\left(\frac{N\tau_1 z}{r}\right)\\
&=\frac{\tau }{2\pi N}\left(T_{1}+T_{2}\right),
\end{align*}
where
\begin{align}
T_{1}&=\int_{0}^{1}\frac{{\rm d}r}{r}\int_{-\infty}^{\infty}{\rm d}z~\psi_{\tau_1,\tau_2}(z)\exp\left[-\tau z^2\left(\frac{1-\frac{r}{\tau_{1}N^2}}{{r}/{\tau_1}}
\right)\right]{\hat f}\left(\frac{N\tau_1 z}{r}\right),\\
T_{2}&=\int_{1}^{N^{2}}\frac{{\rm d}r}{r}\int_{-\infty}^{\infty}{\rm d}z~\psi_{\tau_1,\tau_2}(z)\exp\left[-\tau z^2\left(\frac{1-\frac{r}{\tau_{1}N^2}}{{r}/{\tau_1}}
\right)\right]{\hat f}\left(\frac{N\tau_1 z}{r}\right).\label{3.17'}
\end{align}
We can see that, for all $N\geq1$
\begin{align*}
\int_{0}^{1}\frac{{\rm d}r}{r}\exp\left[-\tau \tau_1 z^2\left(\frac{1-\frac{r}{\tau_{1}N^2}}{{r}}\right)\right]
\leq\int_{0}^{1}\frac{{\rm d}r}{r}\exp\left[-\tau \tau_1 z^2\left(\frac{1-r}{{r}}\right)\right]\leq\log_{+}\left(\frac{e}{\tau \tau_1 z^2}\right),
\end{align*}
where $\log_{+}(x)=\log(e+x)$. Hence,
\begin{align*}
T_{1}\leq f(\mathbb R)\int_{-\infty}^{\infty}\vert\psi_{\tau_1,\tau_2}(z)\vert\log_{+}\left(\frac{e}{\tau \tau_1 z^2}\right){\rm d}z<\infty,
\end{align*}
which implies that
\begin{align}
\limsup_{N\rightarrow\infty}\frac{N}{\log N}\frac{\tau }{2\pi N}T_1=0.\label{3.14}
\end{align}
Then, it remains to prove
$$\liminf_{N\rightarrow\infty}\frac{N}{\log N}\frac{\tau }{2\pi N}T_2\geq t_1f(\mathbb R).$$
We make the change of variable $s={r}{N^{-2}}$ to see that
\begin{align}
T_2
&\geq \int_{-\infty}^{\infty}{\rm d}z\int_{{1}/{N}}^{1}\frac{{\rm d}s}{s}~\psi_{\tau_1,\tau_2}(z)\exp\left[-\tau z^2\left(\frac{1-\frac{s}{\tau_1}}{{sN^2}/{\tau_1}}\right)\right]{\hat f}\left(\frac{\tau_1 z}{sN}\right)\label{3.15'}\\
&=T_{2,1}+T_{2,2},\label{3.15}
\end{align}
where
\begin{align*}
T_{2,1}
&=\int_{-\log N}^{\log N}{\rm d}z\int_{{1}/{N}}^{1}\frac{{\rm d}s}{s}~\psi_{\tau_1,\tau_2}(z)\exp\left[-\tau z^2\left(\frac{1-\frac{s}{\tau_1}}{{sN^2}/{\tau_1}}\right)\right]{\hat f}\left(\frac{\tau_1 z}{sN}\right),\\
T_{2,2}
&=\int_{\vert z \vert>\log N}{\rm d}z\int_{{1}/{N}}^{1}\frac{{\rm d}s}{s}~\psi_{\tau_1,\tau_2}(z)\exp\left[-\tau z^2\left(\frac{1-\frac{s}{\tau_1}}{{sN^2}/{\tau_1}}\right)\right]{\hat f}\left(\frac{\tau_1 z}{sN}\right).
\end{align*}
Because
\begin{align}
\vert T_{2,2} \vert \leq f(\mathbb R)\log N\int_{\vert z \vert>\log N}\vert \psi_{\tau_1,\tau_2}(z)\vert{\rm d}z=o(\log N),\label{3.16}
\end{align}
we focus on the behavior of $T_{2,1}$. First, we can write
\begin{align*}
\psi_{\tau_1,\tau_2}={\rm Re}\{\psi_{\tau_1,\tau_2}\}
+i~{\rm Im}\{\psi_{\tau_1,\tau_2}\},
\end{align*}
where ${\rm Re}\{\psi_{\tau_1,\tau_2}\}$ and ${\rm Im}\{\psi_{\tau_1,\tau_2}\}$ denote the real and imaginary parts of the function $\psi_{\tau_1,\tau_2}$ respectively.
Moreover, we decompose ${\rm Re}\{\psi_{\tau_1,\tau_2}\}$ into two parts such that
\begin{align}
{\rm Re}\{\psi_{\tau_1,\tau_2}(z)\}=
\psi^{+}_{\tau_1,\tau_2}(z)-\psi^{-}_{\tau_1,\tau_2}(z)
\label{3.16'},
\end{align}
where $\psi^{+}_{\tau_1,\tau_2}$ and $\psi^{-}_{\tau_1,\tau_2}$ denote respectively the positive and negative parts of the function $\psi_{\tau_1,\tau_2}$.
A few lines of computation
show that ${\rm Im}\{\psi_{\tau_1,\tau_2}\}$ is an odd function, then by
making a change of variable $r=sN$, we have
\begin{align*}
T_{2,1}
&=\int_{-\log N}^{\log N}{\rm d}z\int_{1}^{N}\frac{{\rm d}r}{r}~\psi_{\tau_1,\tau_2}(z)\exp\left[-\frac{\tau z^2}{N}\left(\frac{\tau_1}{r}-\frac{1}{N}\right)\right]{\hat f}\left(\frac{\tau_1 z}{r}\right)\\
&=T_{2,1,1}-T_{2,1,2},
\end{align*}
where
\begin{align*}
T_{2,1,1}&=\int_{-\log N}^{\log N}{\rm d}z\int_{1}^{N}\frac{{\rm d}r}{r}~\psi^{+}_{\tau_1,\tau_2}(z)\exp\left[-\frac{\tau z^2}{N}\left(\frac{\tau_1}{r}-\frac{1}{N}\right)\right]{\hat f}\left(\frac{\tau_1 z}{r}\right),\\
T_{2,1,2}&=\int_{-\log N}^{\log N}{\rm d}z\int_{1}^{N}\frac{{\rm d}r}{r}~\psi^{-}_{\tau_1,\tau_2}(z)\exp\left[-\frac{\tau z^2}{N}\left(\frac{\tau_1}{r}-\frac{1}{N}\right)\right]{\hat f}\left(\frac{\tau_1 z}{r}\right).
\end{align*}
Notice that
\begin{align*}
T_{2,1,1}\geq\exp\left[-\frac{\tau\tau_1}{N}(\log N)^{2}\right]\int_{-\log N}^{\log N}{\rm d}z\int_{1}^{N}\frac{{\rm d}r}{r}~\psi^{+}_{\tau_1,\tau_2}(z){\hat f}\left(\frac{\tau_1 z}{r}\right).
\end{align*}
Hence,
\begin{align*}
T_{2,1,1}
&\geq(1+o(1))\int_{-\log N}^{\log N}{\rm d}z\int_{(\log N)^2}^{N}\frac{{\rm d}r}{r}~\psi^{+}_{\tau_1,\tau_2}(z){\hat f}\left(\frac{\tau_1 z}{r}\right)\\
&=(f(\mathbb R)+o(1))\int_{-\log N}^{\log N}{\rm d}z\int_{(\log N)^2}^{N}\frac{{\rm d}r}{r}~\psi^{+}_{\tau_1,\tau_2}(z)\\
&=(f(\mathbb R)+o(1))\log N\int_{-\infty}^{\infty}{\rm d}z~\psi^{+}_{\tau_1,\tau_2}(z).
\end{align*}
Recall the definition of $T_{2,1,2}$. We deduce that
\begin{align*}
T_{2,1,2}
&\leq f(\mathbb R)\int_{-\infty}^{\infty}\psi^{-}_{\tau_1,\tau_2}(z){\rm d}z\int_{1}^{N}\frac{{\rm d}r}{r}
=f(\mathbb R)\log N\int_{-\infty}^{\infty}\psi^{-}_{\tau_1,\tau_2}(z){\rm d}z.
\end{align*}
Therefore,
\begin{align}
T_{2,1}\geq f(\mathbb R)\log N\int_{-\infty}^{\infty}\psi_{\tau_1,\tau_2}(z){\rm d}z+o(\log N)=\frac{2\pi f(\mathbb R)}{\tau_2}\log N+o(\log N),\label{3.17}
\end{align}
where we apply the the Parseval's identity in the equality. Combining (\ref{3.14})-(\ref{3.17}), we obtain the first inequality in (\ref{3.15''}).
\end{proof}

\begin{proof}[Proof of Proposition \ref{prop3.2}]
According to Lemma \ref{lem3.6}, it remains to prove the case $\lim_{x\rightarrow\infty}{\hat f}(x)=0$ by showing that
\begin{align*}
\limsup_{N\rightarrow\infty}\frac{N}{\log N}\frac{\tau }{2\pi N}T_2\leq t_{1}f(\mathbb R).
\end{align*}
Recall the definition of $T_2$ from (\ref{3.17'}) that
\begin{align*}
T_2&=\int_{-\infty}^{\infty}{\rm d}z~\psi_{\tau_1,\tau_2}(z)\int_{{1}/{N^2}}^{1}\frac{{\rm d}s}{s}\exp\left[-\tau z^2\left(\frac{1-\frac{s}{\tau_1}}{{sN^2}/{\tau_1}}\right)\right]{\hat f}\left(\frac{\tau_1 z}{sN}\right)\\
&=T'_{2,1}+T'_{2,2},
\end{align*}
where
\begin{align*}
T'_{2,1}&=\int_{-\log N}^{\log N}{\rm d}z~\psi_{\tau_1,\tau_2}(z)\int_{{1}/{N^2}}^{1}\frac{{\rm d}s}{s}\exp\left[-\tau z^2\left(\frac{1-\frac{s}{\tau_1}}{{sN^2}/{\tau_1}}\right)\right]{\hat f}\left(\frac{\tau_1 z}{sN}\right),\\
T'_{2,2}&=\int_{\vert z \vert>\log N}{\rm d}z~\psi_{\tau_1,\tau_2}(z)\int_{{1}/{N^2}}^{1}\frac{{\rm d}s}{s}\exp\left[-\tau z^2\left(\frac{1-\frac{s}{\tau_1}}{{sN^2}/{\tau_1}}\right)\right]{\hat f}\left(\frac{\tau_1 z}{sN}\right).
\end{align*}
Using the same computation as the proof of (\ref{3.16}), it is possible to prove that $T'_{2,2}=o(\log N)$. Then we move on to study the behavior of $T'_{2,1}$. By a change of variable $r=sN$,
\begin{align*}
T'_{2,1}
&=\int_{-\log N}^{\log N}{\rm d}z~\psi_{\tau_1,\tau_2}(z)\int_{{1}/{N}}^{(\log N)^2}\frac{{\rm d}r}{r}\exp\left[-\frac{\tau z^2}{N}\left(\frac{\tau_1}{r}-\frac{1}{N}\right)\right]{\hat f}\left(\frac{\tau_1 z}{r}\right)\\
&~~~+\int_{-\log N}^{\log N}{\rm d}z~\psi_{\tau_1,\tau_2}(z)\int_{(\log N)^2}^{N}\frac{{\rm d}r}{r}\exp\left[-\frac{\tau z^2}{N}\left(\frac{\tau_1}{r}-\frac{1}{N}\right)\right]{\hat f}\left(\frac{\tau_1 z}{r}\right)\\
&:=T'_{2,1,1}+T'_{2,1,2}.
\end{align*}
According to the proof of Theorem 5.2, item 2 in \cite{MR4242879}, since $\lim_{x\rightarrow\infty}{\hat f}(x)=0$ and $\int_{-\infty}^{\infty}\vert \psi_{\tau_1,\tau_2}(z)\vert{\rm d}z<\infty$, it is easy to see that $T'_{2,1,1}=o(\log N)$. Recall the function $\psi^{+}_{\tau_1,\tau_2}$ and $\psi^{-}_{\tau_1,\tau_2}$ defined in (\ref{3.16'}). In order to estimate $T'_{2,1,2}$, we rewrite one more time to find that
\begin{align*}
T'_{2,1,2}
&=\int_{-\log N}^{\log N}{\rm d}z~\psi^{+}_{\tau_1,\tau_2}(z)\int_{(\log N)^2}^{N}\frac{{\rm d}r}{r}\exp\left[-\frac{\tau z^2}{N}\left(\frac{\tau_1}{r}-\frac{1}{N}\right)\right]{\hat f}\left(\frac{\tau_1 z}{r}\right)\\
&~~~-\int_{-\log N}^{\log N}{\rm d}z~\psi^{-}_{\tau_1,\tau_2}(z)\int_{(\log N)^2}^{N}\frac{{\rm d}r}{r}\exp\left[-\frac{\tau z^2}{N}\left(\frac{\tau_1}{r}-\frac{1}{N}\right)\right]{\hat f}\left(\frac{\tau_1 z}{r}\right).
\end{align*}
Then applying a similar argument in proving (\ref{3.17}), which we skip, we conclude that
\begin{align*}
\limsup_{N\rightarrow\infty}\frac{N}{\log N}\frac{\tau }{2\pi N}T'_{2,1,2}\leq t_{1}f(\mathbb R).
\end{align*}
This proves Proposition \ref{prop3.2}.
\end{proof}
\begin{remark}
{\rm Comparing with the proof of asymptotic variance {\rm{\cite{MR4242879}}}, one of the difficulties of our estimation is that $\psi_{\tau_1,\tau_2}$ may not be a real function and the real part of the function $\psi_{\tau_1,\tau_2}$ might be negative unless $\tau_1=\tau_2$. So we use highly technical computations to obtain the results. For example, the values of the expressions can be reduced or enlarged by adjusting the upper and lower limits of integration with respect to time variable, rather than the space variable; see $(\ref{3.16})$. And we need to decompose the function $\psi_{\tau_1,\tau_2}$ in two parts; see $(\ref{3.16'})$.}
\end{remark}

\subsection{Analysis of Riesz kernel case}
\begin{proof}[Proof of Proposition \ref{prop3.3}]{\rm\bf Case 1:} $\beta\in(0,1).$
Recall (\ref{3.9}). Since $f(dx)=\Vert x \Vert^{-\beta}$ we can write
\begin{align}
V^{(1)}_{N,t_1,t_2}
&=\int_{{{\mathbb R}^{d}}}{\rm d}z\left(I_{\tau_2}\ast {\tilde I_{\tau_1}}\right)(z)\int_{0}^{t_{1}}{\rm d}s~\E\left(\left\Vert \sqrt{\frac{2s(\tau-s)}{\tau}}Z+\frac{Ns}{\tau}z
\right\Vert^{-\beta}\right)\notag\\
&=\frac{1}{N^\beta}\int_{{{\mathbb R}^{d}}}{\rm d}z\left(I_{\tau_2}\ast {\tilde I_{\tau_1}}\right)(z)\int_{0}^{t_{1}}{\rm d}s~\E\left(\left\Vert \frac{1}{N}\sqrt{\frac{2s(\tau-s)}{\tau}}Z+\frac{s}{\tau}z
\right\Vert^{-\beta}\right),\label{3.18}
\end{align}
where $Z$ is a $d$-dimensional standard normal random variable.
Since $\int_{{{\mathbb R}^{d}}}\Vert z\Vert^{-\beta}\left(I_{\tau_2}\ast {\tilde I_{\tau_1}}\right)(z)~{\rm d}z<\infty$, by using Lemma 3.1 in \cite{MR4098872} and the dominated convergence theorem, we have
\begin{align*}
\lim_{N\rightarrow\infty}N^\beta V^{(1)}_{N,t_1,t_2}&=\int_{{{\mathbb R}^{d}}}{\rm d}z\left(I_{\tau_2}\ast {\tilde I_{\tau_1}}\right)(z)\int_{0}^{t_{1}}{\rm d}s~\tau^\beta s^{-\beta} \Vert z\Vert^{-\beta}.
\end{align*}
Next, recall (\ref{3.11}). Before deducing the last two parts, we rewrite the expression $V^{(1)}_{N,t_1,t_2}$  to see that
\begin{align}
&V^{(1)}_{N,t_1,t_2}
=\frac{\kappa_{\beta,d}}{(2\pi)^d}\int_{0}^{t_{1}}
{{\rm d}s}\int_{{\mathbb R}^d}{\rm d}z~\psi_{\tau_1,\tau_2}(z)
\left(\frac{\tau}{sN}\right)^d\exp\left(-\frac{\tau(\tau-s)}{N^2s}\Vert z\Vert^2\right)\left\Vert\frac{\tau}{sN}z
\right\Vert^{\beta-d}\notag\\
&=\frac{\kappa_{\beta,d}}{(2\pi)^d N^\beta}\int_{0}^{t_{1}}\left(\frac{\tau}{s}\right)^\beta{{\rm d}s}
\int_{{\mathbb R}^d}{\rm d}z~\psi_{\tau_1,\tau_2}(z)\exp\left(-\frac{\tau(\tau-s)}{N^2s}\Vert z\Vert^2\right)\left\Vert z\right\Vert^{\beta-d} \label{3.20}\\
&=\frac{\tau\kappa_{\beta,d}}{(2\pi)^d N^\beta}\int_{{\mathbb R}^d}{\rm d}z\left\Vert z\right\Vert^{\beta-d}\psi_{\tau_1,\tau_2}(z)
\int_{\tau_1-1}^{\infty}{{\rm d}r}~(1+r)^{\beta-2}\exp\left(-\frac{r\tau}{N^2}\Vert z\Vert^2\right) \notag\\
&=\frac{\tau^{2-\beta}\kappa_{\beta,d}}{(2\pi)^d N^{2-\beta}}\int_{{\mathbb R}^d}{\rm d}z\left\Vert z\right\Vert^{2-\beta-d}\psi_{\tau_1,\tau_2}(z)
\int_{\frac{\left({\tau_1-1}\right)\tau{\Vert z \Vert}^2}{N^2}}^{\infty}{{\rm d}r}\left(\frac{\tau\Vert z \Vert^2}{N^2}+r\right)^{\beta-2}e^{-r},\label{3.21}
\end{align}
where we use scaling property of Fourier transform in the first equality and change of variables in the last two equalities.
Now, we move on to prove the last two parts. \\
{\rm\bf Case 2:} $\beta=1$.
We apply (\ref{3.20}) with $\beta=1$ to see that
\begin{align*}
\frac{N}{\log N}V^{(1)}_{N,t_1,t_2}
&=\frac{\kappa_{1,d}}{(2\pi)^d \log N}\int_{0}^{t_{1}}\frac{\tau}{s}{{\rm d}s}
\int_{{\mathbb R}^d}{\rm d}z~\psi_{\tau_1,\tau_2}(z)\exp\left(-\frac{\tau(\tau-s)}{N^2s}\Vert z\Vert^2\right)\left\Vert z\right\Vert^{1-d}\\
&=\frac{\tau\kappa_{1,d}}{(2\pi)^d\log N}\int_{0}^{\tau}\frac{{{\rm d}s}}{s}
\int_{{\mathbb R}^d}{\rm d}z~\psi_{\tau_1,\tau_2}(z)\exp\left(-\frac{\tau(\tau-s)}{N^2s}\Vert z\Vert^2\right)\left\Vert z\right\Vert^{1-d}\\
&~~~-\frac{\tau\kappa_{1,d}}{(2\pi)^d\log N}\int_{t_1}^{\tau}\frac{{{\rm d}s}}{s}
\int_{{\mathbb R}^d}{\rm d}z~\psi_{\tau_1,\tau_2}(z)\exp\left(-\frac{\tau(\tau-s)}{N^2s}\Vert z\Vert^2\right)\left\Vert z\right\Vert^{1-d}\\
&:=\frac{1}{\log N}(A_1-A_2).
\end{align*}
According to Lemma \ref{lem6.4} and the dominated convergence theorem, the following:
\begin{align}
\lim_{N\rightarrow\infty}\frac{A_1}{\log N}=\frac{2\tau\kappa_{1,d}}{(2\pi)^d}\int_{{\mathbb R}^d}\psi_{\tau_1,\tau_2}(z)\left\Vert z\right\Vert^{1-d}{\rm d}z\label{3.24'}
\end{align}
holds due to $\int_{{\mathbb R}^d}\left\Vert z\right\Vert^{1-d}\vert\psi_{\tau_1,\tau_2}(z)\vert
\log_+(1/\Vert z\Vert){\rm d}z<\infty.$
Moreover, we observe that
\begin{align*}
0\leq A_2\leq\frac{\tau\kappa_{1,d}}{(2\pi)^d}\log {\left(\frac{\tau}{t_1}\right)}\int_{{\mathbb R}^d}\vert \psi_{\tau_1,\tau_2}(z)\vert\left\Vert z\right\Vert^{1-d}{\rm d}z=o(\log N),
\end{align*}
which together with (\ref{3.24'}) proves (\ref{3.6}).\\
{\rm\bf Case 3:} $\beta\in(1,2\wedge d)$. Recall (\ref{3.21}). Since $\int_{{\mathbb R}^d}\left\Vert z\right\Vert^{2-\beta-d}\vert \psi_{\tau_1,\tau_2}(z)\vert {\rm d}z<\infty$ and $\int_{0}^{\infty}{r}^{\beta-2}e^{-r}<\infty$, the dominated convergence theorem implies that
\begin{align*}
\lim_{N\rightarrow\infty}N^{2-\beta}V^{(1)}_{N,t_1,t_2}=
\frac{{\tau^{2-\beta}\kappa_{\beta,d}}}{{2\pi}^d}\int_{{\mathbb R}^d}{{\psi_{\tau_1,\tau_2}}}(z){\Vert z \Vert}^{2-\beta-d}{\rm d}z\int_{0}^\infty{r^{\beta-2}}e^{-r}{\rm d}r.
\end{align*}
The proof is completed.
\end{proof}

\section{Tightness}
We will show the upper bound and tightness of the spatial averages in this section. To simplify the expressions, let $L$ denote a real number which may take different values and depend on different variables, this notation will also be used in Chapter 5.
\subsection{Uniform upper bound}
In this subsection, we prove quantitative upper bound of the spatial averages in $L^k(\Omega)$. Apparently, the rates of convergence are coincident with the rates that were ensured by Propositions \ref{prop3.1}-\ref{prop3.3}.
\begin{lemma}\label{lem4.1}
${(d\geq2)}$.
 Suppose ${\mathcal R}(f)<\infty$, then for any $k\geq2$ and $T\geq0$,
\begin{align}
\sup_{N\geq e}\left\Vert \sqrt{N} S_{N,t}\right\Vert_k\leq L_{T,k}\sqrt{t}~~~{\it uniformly~for~all}~t\in(0,T).\label{4.1}
\end{align}
\end{lemma}
\begin{lemma}\label{lem4.2}
{$(d=1)$}. Suppose $f(\mathbb R)<\infty$, then for any $k\geq2$ and $T\geq0$,
\begin{align}
\sup_{N\geq e}\left\Vert \sqrt{\frac{N}{\log N}} S_{N,t}\right\Vert_k \leq L_{T,k}\sqrt{t\log_{+}{\left(1/t\right)}}~~~{\it uniformly~for~all}~t\in(0,T).\label{4.2}
\end{align}
\end{lemma}
\begin{lemma}\label{lem4.3}
{\rm(Riesz kernel)}. Suppose $f({\rm d}x)=\Vert x\Vert^{-\beta}$, where $0<\beta<2\wedge d$, then for any $k\geq2$ and $T\geq0$, we have\\
$\mathbf{(A)}$ If $\beta\in(0,1)$, then
\begin{align}
\sup_{N\geq e}\left\Vert \sqrt{N^\beta} S_{N,t}\right\Vert_k \leq L_{T,k,\beta}\sqrt{t}~~~{\it uniformly~for~all}~t\in(0,T),\label{4.3}
\end{align}
$\mathbf{(B)}$ If $\beta=1$, then
\begin{align}
\sup_{N\geq e}\left\Vert \sqrt{\frac{N}{\log N}} S_{N,t}\right\Vert_k \leq L_{T,k}\sqrt{t}~~~{\it uniformly~for~all}~t\in(0,T),\label{4.4}
\end{align}
$\mathbf{(C)}$ If $\beta\in(1,2\wedge d)$, then
\begin{align}
\sup_{N\geq e}\left\Vert \sqrt{N^{2-\beta}} S_{N,t}\right\Vert_k \leq L_{T,k,\beta}\sqrt{t^{2-\beta}}~~~{\it uniformly~for~all}~t\in(0,T).\label{4.5}
\end{align}
\end{lemma}
\begin{proof}[Proof of Lemmas \ref{lem4.1}-\ref{lem4.3}]
Denote $V^{(1)}_{N,t}:=V^{(1)}_{N,t,t}$, and define
\begin{align}
g_{N,t}(s,y):=\frac{1}{N^{d}}\int_{[0,N]^d}{\boldsymbol p}_{s(t-s)/t}{\left(y-\frac{s}{t}x\right){\rm d}x}.\label{4.6}
\end{align}
Recall (\ref{1.4}) and (\ref{1.5}). Thanks to stochastic Fubini's theorem and the BDG inequality,
\begin{align*}
\Vert S_{N,t}\Vert_k^2&=\left\Vert \int_{(0,t)\times{\mathbb R}^d}\left(\frac{1}{N^{d}}\int_{[0,N]^d}{\boldsymbol p}_{s(t-s)/t}{\left(y-\frac{s}{t}x\right){\rm d}x}\right)U(s,y)\eta(ds,dy)\right\Vert_k^2\\
&\leq{L_k}\int_{(0,t)\times{\mathbb R}^{2d}}{{\rm d}y}f({{\rm d}y'}){{\rm d}s}~g_{N,t}(s,y)g_{N,t}(s,y+y')\Vert U(s,y)U(s,y+y')\Vert_{k/2}.
\end{align*}
Then, we use the moment inequality (\ref{2.2}) and the semigroup property of the heat kernel to find that
\begin{align*}
\Vert S_{N,t}\Vert_k^2
&\leq\frac{L_{T,k}}{N^{2d}}\int_{(0,t)\times{\mathbb R}^{d}}f({{\rm d}y'}){{\rm d}s}\int_{[0,N]^{2d}}{{\rm d}x_1}{{\rm d}x_2}~{\boldsymbol p}_{2s(t-s)/t}\left(y'-\frac{s}{t}(x_2-x_1)\right)\\
&\leq{}{L_{T,k}}\left\vert V_{N,t}^{(1)}\right\vert.
\end{align*}
First, we analyze the case that $d\geq2$. According to the proof of proposition \ref{prop3.1}, we observe that
\begin{align*}
\Vert S_{N,t}\Vert_k^2
&\leq L_{T,k}\frac{t}{N}\int_0^\infty \frac{{{\rm d}s}}{s^d}\int_{{\mathbb R}^{d}}{\hat f}({{\rm d}z})\left\vert\widehat{\left[\left(I_{1}\ast {\tilde I_{1}}\right)\left(\frac{\bullet}{s}\right)\right]}(z)\right\vert\\
&=L_{T,k}\frac{t}{N}.
\end{align*}
This proves Lemma \ref{lem4.1}, and the rest Lemma \ref{lem4.2} and Lemma \ref{lem4.3} can be proved by using similar arguments, we skip the details.
\end{proof}
\subsection{Moment estimate}
In this subsection, we establish tightness by the following moment estimates.
\begin{proposition}\label{prop4.5}
${(d\geq2)}$. Assume ${\mathcal R}(f)<\infty$, then for any $T>0$, $k\geq2$, and $\gamma\in(0,{1}/{8})$, there exists a real number $L=L_{T,k,\gamma}>0$ such that for all $\varepsilon\in(0,1)$,
\begin{align}
\sup_{t\in(0,T)}\E\left(\left\vert S_{N,t+\varepsilon}-S_{N,t}\right\vert^k\right)\leq L\varepsilon^{\gamma k}N^{-k/2}~~~{\it uniformly~for~all~}N\geq e.\label{4.7}
\end{align}
\end{proposition}
\begin{proposition}\label{prop4.6}
${(d=1)}$. Assume $f(\mathbb R)<\infty$, then for any $T>0$, $k\geq2$, and $\gamma\in(0,{1}/{4})$, there exists a real number $L=L_{T,k,\gamma}>0$ such that for all $\varepsilon\in(0,1)$,
\begin{align}
\sup_{t\in(0,T)}\E\left(\left\vert S_{N,t+\varepsilon}-S_{N,t}\right\vert^k\right)\leq L\varepsilon^{\gamma k}\left(\frac{N}{\log N}\right)^{-k/2}~~~{\it uniformly~for~all~}N\geq e.\label{4.8}
\end{align}
\end{proposition}
\begin{proposition}\label{prop4.7}
${\rm(Riesz~kernel)}$. Assume $f({\rm d}x)=\Vert x \Vert^{-\beta}{\rm d}x$, $\beta\in (0,2\wedge d)$, for any $T>0$, $k\geq2$ and $\varepsilon\in(0,1)$ we have\\
$\mathbf{(A)}$ If $\beta\in(0,1)$, then for any $\gamma\in(0,1/4]$, there exists $L=L_{T,k,\gamma,\beta}>0$ such that
\begin{align}
\sup_{t\in(0,T)}\E\left(\left\vert S_{N,t+\varepsilon}-S_{N,t}\right\vert^k\right)\leq L\varepsilon^{\gamma k}\left(N\right)^{-\beta k/2}~~~{\it uniformly~for~all~}N\geq e.\label{4.9}
\end{align}
$\mathbf{(B)}$ If $\beta=1$, then for any $\gamma\in(0,1/4)$, there exists $L=L_{T,k,\gamma}>0$ such that
\begin{align}
\sup_{t\in(0,T)}\E\left(\left\vert S_{N,t+\varepsilon}-S_{N,t}\right\vert^k\right)\leq L\varepsilon^{\gamma k}\left(\frac{N}{\log N}\right)^{- k/2}~~~{\it uniformly~for~all~}N\geq e.\label{4.10}
\end{align}
$\mathbf{(C)}$ If $\beta\in(1,2\wedge d)$, then for any $\gamma\in(0, (2-\b)/(6-2\b)]$, there exists $L=L_{T,k,\gamma,\beta}>0$ such that
\begin{align}
\sup_{t\in(0,T)}\E\left(\left\vert S_{N,t+\varepsilon}-S_{N,t}\right\vert^k\right)\leq L\varepsilon^{\gamma k}\left(N\right)^{-(2-\beta)k/2}~~~{\it uniformly~for~all~}N\geq e.\label{4.11}
\end{align}
\end{proposition}
The proofs of Propositions \ref{prop4.5}-\ref{prop4.7} hinge on the following lemmas, which will be used to analyze the behavior when $t$ stays away from zero.
\begin{lemma}\label{lem4.8}
$(d\geq2)$. Suppose ${\mathcal R}(f)<\infty$, then for any $T>0$, $k\geq2$, and $\alpha\in(0,{1}/{4})$, there exists a real number $L=L_{T,k,\alpha}>0$ such that for all $\varepsilon\in(0,1)$, $t\in(0,T)$ and $N\geq e$,
\begin{align}
\E\left(\left\vert S_{N,t+\varepsilon}-S_{N,t}\right\vert^k\right)\leq \frac{L\varepsilon^{\alpha k}}{(t\wedge1)^{k/2}}N^{-k/2}.\label{4.12}
\end{align}
\end{lemma}
\begin{lemma}\label{lem4.99}
Suppose $f(\R)<\infty$, then for any $T>0$, $k\geq2$, and $\delta>0$, there exists a real number $L=L_{T,k,\delta}>0$ such that for all  $\varepsilon\in(0,1)$, $t\in(0,T)$ and $N\geq e$,
\begin{align}
\E\left(\left\vert S_{N,t+\varepsilon}-S_{N,t}\right\vert^k\right)\leq \frac{L\varepsilon^{k/2}}{(t\wedge1)^{k(1+\delta)/2}}\left(\frac{N}{\log N}\right)^{-k/2}.\label{4.14}
\end{align}
\end{lemma}
\begin{lemma}\label{lem4.9}
${\rm(Riesz~kernel)}$. Suppose $f({\rm d}x)=\Vert x \Vert^{-\beta}{\rm d}x$, $\beta\in (0,2\wedge d)$, for any $T>0$, $k\geq2$ and $\varepsilon\in(0,1)$ we have\\
$\mathbf{(A)}$ If $\beta\in(0,1)$, then there exists a real number $L=L_{T,k,\beta}>0$ such that for all $t\in(0,T)$ and $N\geq e$,
\begin{align}
\E\left(\left\vert S_{N,t+\varepsilon}-S_{N,t}\right\vert^k\right)\leq
\frac{L\varepsilon^{k/2}}{(t\wedge1)^{k/2}}N^{-\beta k/2}.\label{4.13}
\end{align}
$\mathbf{(B)}$ If $\beta=1$, then for any $\delta>0$, there exists a real number $L=L_{T,k,\delta}>0$ such that for all $t\in(0,T)$ and $N\geq e$,
\begin{align}
\E\left(\left\vert S_{N,t+\varepsilon}-S_{N,t}\right\vert^k\right)\leq \frac{L\varepsilon^{k/2}}{(t\wedge1)^{k(1+\delta)/2}}\left(\frac{N}{\log N}\right)^{-k/2}.\label{4.14}
\end{align}
$\mathbf{(C)}$ If $\beta\in(1,2\wedge d)$, then there exists a real number $L=L_{T,k,\beta}>0$ such that for all $t\in(0,T)$ and $N\geq e$,
\begin{align}
\E\left(\left\vert S_{N,t+\varepsilon}-S_{N,t}\right\vert^k\right)\leq \frac{L\varepsilon^{ k/2}}{(t\wedge1)^{k/2}}N^{-(2-\beta)k/2}.
\label{4.15}
\end{align}
\end{lemma}
We now aim to prove Lemmas \ref{lem4.8}-\ref{lem4.9}. Recall from (\ref{1.4}) and (\ref{1.5}) that
\begin{align*}
S_{N,t+\varepsilon}-&S_{N,t}
=\frac{1}{N^d}\int_{[0,N]^d}[U(t+\varepsilon,x)-U(t,x)]{\rm d}x\\
=&\int_{(0,t)\times \mathbb R^d}U(s,y){\mathcal A}(s,y)\eta({\rm d}s,{\rm d}y)+\int_{(t,t+\varepsilon)\times \R^d}U(s,y){\mathcal B}(s,y)\eta({\rm d}s,{\rm d}y),
\end{align*}
where
\begin{align*}
{\mathcal A}(s,y)&=\frac{1}{N^d}\int_{[0,N]^d}\left[{\boldsymbol p}_{s(t+\varepsilon-s)/t+\varepsilon}
\left(y-\frac{sx}{t+\varepsilon}\right)
-{\boldsymbol p}_{s(t-s)/t}\left(y-\frac{sx}{t}\right)\right]{\rm d}x,\\
{\mathcal B}(s,y)&=\frac{1}{N^d}\int_{[0,N]^d}{\boldsymbol p}_{s(t+\e-s)/(t+\e)}
\left(y-\frac{sx}{t+\e}\right){\rm d}x.
\end{align*}
Thus,
\begin{align}
\Vert S_{N,t+\varepsilon}-S_{N,t}\Vert_k\leq T_{{\mathcal A}}+T_{{\mathcal B}},\label{4.16}
\end{align}
where
\begin{align*}
T_{{\mathcal A}}&=\left\Vert\int_{(0,t)\times \R^d}U(s,y){\mathcal A}(s,y)\eta({\rm d}s,{\rm d}y)\right\Vert_k\\
~{\rm and}~T_{{\mathcal B}}&=\left\Vert\int_{(t,t+\varepsilon)\times \R^d}U(s,y){\mathcal B}(s,y)\eta({\rm d}s,{\rm d}y)\right\Vert_k.
\end{align*}
Owing to the BDG inequality, (\ref{2.2}) and Parseval's identity,
\begin{align}
T_{{\mathcal A}}^2&\leq\frac{L_{T,k}}{(2\pi)^d}\int_0^t{\rm d}s\int_{{\mathbb R}^d}{\hat f}({\rm d}\xi)\left \vert{\widehat{{\mathcal A}(s)}(\xi)}\right\vert^2,\label{4.17}\\
T_{{\mathcal B}}^2&\leq L_{T,k}\int_{(t,t+\varepsilon)\times \R^{2d}}{\rm d}yf({\rm d}y'){\rm d}s~{\mathcal B}(s,y){\mathcal B}(s,y+y').\label{4.18}
\end{align}
\begin{proof}[Proof of Lemma \ref{lem4.8}] We will estimate $T_{{\mathcal A}}$ and
$T_{{\mathcal B}}$ separately.
First, by using a change of variable [$s\mapsto(t+\e)s/N$] we can see that
\begin{align*}
T_{{\mathcal A}}^2
\leq\frac{L_{T,k}(t+\varepsilon)}{N}
\int_0^{tN/(t+\varepsilon)}{\rm d}s
\int_{{\mathbb R}^d}{\hat f}({\rm d}\xi)
\left\vert{\widehat{{\mathcal A}\left({(t+\varepsilon)s}/{N}\right)}(\xi)}
\right\vert^2.
\end{align*}
Notice that
\begin{align}
{\widehat{{\mathcal A}\left({(t+\varepsilon)s}/{N}\right)}(\xi)}
&=\frac{1}{N^d}\int_{[0,N]^d}\left[
\exp\left(\frac{isx\cdot\xi}{N}
-\frac{s(t+\varepsilon)(N-s)\Vert \xi \Vert^2}{2N^2}\right.
\right)\notag\\
%{\rm d}x\\
&~~~-
%\frac{1}{N^d}\int_{[0,N]^d}
\left.\exp\left(\frac{i(t+\varepsilon)sx\cdot\xi}{Nt}-
\frac{s(t+\varepsilon)(Nt-(t+\e)s)\Vert \xi \Vert^2}{2N^2t}\right)\right]{\rm d}x\notag\\
&=J_1+J_2,\label{4.19}
\end{align}
where
\begin{align*}
J_1&=\int_{[0,1]^d}e^{isx\cdot\xi}{\rm d}x~
\left[\exp\left(-\frac{\varepsilon(t+\varepsilon)s^2\Vert \xi \Vert^2}{2N^2t}\right)-1\right]
\exp\left[-\frac{s(t+\varepsilon)\Vert \xi \Vert^2}{2N}\left(1-\frac{s(t+\varepsilon)}{Nt}\right)\right],\\
J_2&=\int_{[0,1]^d}\left[\exp(isx\cdot\xi)-\exp\left(\frac{is(t+\varepsilon)x\cdot\xi}{t}\right)\right]{\rm d}x~
\exp\left[-\frac{s(t+\varepsilon)\Vert \xi \Vert^2}{2N}\left(1-\frac{s(t+\varepsilon)}{Nt}\right)\right].
\end{align*}
Hence,
\begin{align}
T_{{\mathcal A}}^2
&\leq \frac{L_{T,k}}{N}\int_0^{\infty}{\rm d}s
\int_{{\mathbb R}^d}{\hat f}({\rm d}\xi)
\left\vert J_1+J_2\right\vert^2\notag\\
&\lesssim \frac{L_{T,k}}{N}\int_0^{\infty}{\rm d}s
\int_{{\mathbb R}^d}{\hat f}({\rm d}\xi)
\left\vert J_1\right\vert^2+\frac{L_{T,k}}{N}\int_0^{\infty}{\rm d}s
\int_{{\mathbb R}^d}{\hat f}({\rm d}\xi)
\left\vert J_2\right\vert^2\notag\\
&:=T_{{\mathcal A},1}^2+T_{{\mathcal A},2}^2.\label{4.20}
\end{align}
Next, we estimate $T_{{\mathcal A}}^2$ by analyzing the behavior of $\vert J_1\vert^2$ and $\vert J_2\vert^2$. Define
\begin{align}
\phi(x):=\frac{1-\cos x}{x^2}~~~{\rm for~all~}x\in{\mathbb R}\backslash\{0\},\label{4.30}
\end{align}
and $\phi(0):=1/2$ to preserve continuity. We have,
\begin{align*}
\vert J_1\vert^2
&\leq\left\vert\int_{[0,1]^d}e^{isx\cdot\xi}{\rm d}x\right\vert^2
\left\vert1-\exp\left(-\frac{\varepsilon(t+\varepsilon)s^2\Vert \xi \Vert^2}{2N^2t}\right)\right\vert^2\\
&\lesssim L_{\alpha}\prod_{j=1}^d
\phi(s\xi_j)\frac{[\varepsilon(t+\varepsilon)]^{2\alpha}(s\Vert \xi \Vert)^{4\alpha}}{(N^2t)^{2\alpha}}\lesssim L_{\alpha}\left[\prod_{j=1}^d\left(1\wedge \frac{1}{s\vert \xi_j \vert^2}\right)\right]\frac{\varepsilon^{2\alpha}(s\Vert \xi \Vert)^{4\alpha}}{(N^2t)^{2\alpha}}\\
&\lesssim L_{\alpha}\left(1\wedge{\frac{1}{s\Vert \xi \Vert^2}}\right)\frac{\varepsilon^{2\alpha}(s\Vert \xi \Vert)^{4\alpha}}{(N^2t)^{2\alpha}},
\end{align*}
where in the second inequality we use the fact that for any $\alpha\in(0,1/4)$ and $x>0$, there exists a real number $C_\alpha$  such that
\begin{align}
1-e^{-x}\leq C_\alpha x^\alpha.\label{4.21}
\end{align}
And the last inequality holds by the following:
\begin{align}
\prod_{j=1}^d(1\wedge\vert x_j\vert^{-1})\lesssim1\wedge\Vert x\Vert^{-1}.\label{4.22}
\end{align}
Thus, using a change of variable $[s\mapsto s/\Vert \xi \Vert]$, we have
\begin{align*}
T_{A,1}^2
\lesssim L_{T,k,\alpha}\frac{\varepsilon^{2\alpha}
}{t^{2\alpha}N^{4\alpha+1}}
\int_0^{\infty}{\rm d}s~s^{4\alpha}\left(1\wedge{s}^{-2}\right)
\times\int_{{\mathbb R}^d}\frac{1}{\Vert \xi \Vert}{\hat f}({\rm d}\xi).
%\lesssim L_{T,k,\alpha}\frac{\varepsilon^{2\alpha}(t+\varepsilon)^{2\alpha+1}}{t^{2\alpha}N^{2\alpha+1}},
%\label{4.23}\tag{4.23}
\end{align*}
Since ${\mathcal R}(f)<\infty$, applying Lemma 5.9 in \cite{MR4242879} we see that $\int_{{\mathbb R}^d}{\Vert \xi \Vert}^{-1}{\hat f}({\rm d}\xi)<\infty$, and notice that for every $\alpha\in(0,1/4)$, $\int_0^{\infty}s^{4\alpha}\left(1\wedge s^{-2}\right){\rm d}s<\infty$, we finally get
\begin{align}
T_{A,1}^2\lesssim L_{T,k,\alpha}\frac{\varepsilon^{2\alpha}
}{t^{2\alpha}N^{4\alpha+1}}.
\label{4.23}
\end{align}
Now, we turn to analyze $\vert J_2\vert^2$, a few lines of computation show that
\begin{align*}
\vert J_2\vert^2
&\leq \left\vert\int_{[0,1]^d}\left(e^{i\left(1+\varepsilon/t\right)s x\cdot\xi}-e^{isx\cdot\xi}\right){\rm d}x\right\vert^2\\
&\leq \prod_{j=1}^{d}\left\vert\frac{t}{i(t+\varepsilon)s\xi_j}e^{is\xi_j}\left(e^{i\varepsilon s\xi_j/t}-1\right)+\frac{\varepsilon}{i(t+\varepsilon)s\xi_j}\left(1-e^{is\xi_j}\right)\right\vert^2\\
&\leq
\prod_{j=1}^{d}\left\{\frac{2t^2}{(t+\varepsilon)^2(s\xi_j)^2}\left[2-2\cos(\varepsilon s\xi_j/t)\right]+\frac{2\varepsilon^2}{(t+\varepsilon)^2(s\xi_j)^2}\left[2-2\cos( s\xi_j)\right]\right\}.
\end{align*}
Since $\vert1-\cos x\vert\leq 2\wedge (x^2/2)$, for every $x\in\R$, we rescale one more time to see that
\begin{align*}
\vert J_2\vert^2
&\leq
\prod_{j=1}^{d}\left\{\frac{2t^2}{(t+\varepsilon)^2(s\xi_j)^2}
\left[4\wedge\left(\frac{\varepsilon^2(s\xi_j)^2}{t^2}\right)\right]+\frac{2\varepsilon^2}{(t+\varepsilon)^2(s\xi_j)^2}
\left[4\wedge( s\xi_j)^2\right]\right\}\\
&\leq
\prod_{j=1}^{d}\left[\frac{8(t^2+\varepsilon^2)}{(t+\varepsilon)^2(s\xi_j)^2}
\wedge\frac{4\varepsilon^2}{(t+\varepsilon)^2}\right]
\lesssim
\left[\frac{(t^2+\varepsilon^2)}{(t+\varepsilon)^2(s\Vert\xi\Vert)^2}
\wedge\frac{\varepsilon^2}{(t+\varepsilon)^2}\right]\\
&\lesssim\frac{1}{(s\Vert\xi\Vert)^2}\wedge\frac{\varepsilon}{t}\lesssim\left(1\vee\frac{1}{t}\right)\left(\frac{1}{(s\Vert\xi\Vert)^2}\wedge\varepsilon\right),
\end{align*}
where we use (\ref{4.22}) in the third inequality. Therefore, with a change of variable $[s\mapsto s/\Vert \xi \Vert]$ again, we obtain
\begin{align}
T_{A,2}^2
&\lesssim L_{T,k}\frac{(t\vee 1)}{Nt}
\int_0^{\infty}{\rm d}s~\left({s}^{-2}\wedge\e\right)
\times\int_{{\mathbb R}^d}\frac{1}{\Vert \xi \Vert}{\hat f}({\rm d}\xi)\notag\\
&\lesssim L_{T,k,\alpha}\frac{\sqrt{\varepsilon}(t\vee 1)}{Nt},
\label{4.24}
\end{align}
where the last inequality follows from Lemma \ref{lem6.2} (1). Combining (\ref{4.20}), (\ref{4.23}) and (\ref{4.24}), we conclude that
\begin{align}
T_{A}^2\lesssim L_{T,k,\alpha}\frac{\varepsilon^{2\alpha}}{N(t\wedge1)}.
\label{4.25}
\end{align}
Recall {(\ref{4.18})}. In order to estimate $T_{\mathcal B}^2$, we use the semigroup property of the heat kernel and change of variables as follows:
\begin{align*}
T_{\mathcal B}^2
&\leq\frac{L_{T,k}}{N^{2d}}\int_{t}^{t+\varepsilon}{\rm d}s\int_{{\mathbb R}^d}f({\rm d}y')\int_{[0,N]^{2d}}{\rm d}x_1{\rm d}x_2~ {\boldsymbol p}_{2s(t+\varepsilon-s)/(t+\varepsilon)}\left(y'-\frac{s}{t+\varepsilon}(x_2-x_1)\right)\\
&={L_{T,k}}\int_{t}^{t+\varepsilon}{\rm d}s\int_{{\mathbb R}^d}f({\rm d}y')\int_{[0,1]^{2d}}{\rm d}x_1{\rm d}x_2~ {\boldsymbol p}_{2s(t+\varepsilon-s)/(t+\varepsilon)}\left(y'-\frac{sN}{t+\varepsilon}(x_2-x_1)\right)\\
&={L_{T,k}}\frac{t+\varepsilon}{N}\int_{Nt/(t+\varepsilon)}^{N}{\rm d}s\int_{{\mathbb R}^d}f({\rm d}y')\int_{[0,1]^{2d}}{\rm d}x_1{\rm d}x_2~ {\boldsymbol p}_{\gamma,N}\left(y'-s(x_2-x_1)\right),
\end{align*}
where ${{\gamma,N=[2s(t+\varepsilon)(1-s/N)]/N}}$. Furthermore, thanks to Parseval's identity and (\ref{4.22}), we may write the following:
\begin{align}
T_{\mathcal B}^2
&\leq\frac{L_{T,k}}{(2\pi)^dN}\int_{Nt/(t+\varepsilon)}^{N}{\rm d}s\int_{{\mathbb R}^d}{\hat f}({\rm d}\xi)~e^{-\frac{\gamma,N}{2}\Vert \xi\Vert^2}\int_{[0,1]^{2d}}{\rm d}x_1{\rm d}x_2~\exp\left(is(x_2-x_1)\cdot\xi\right)\notag\\
&\lesssim L_{T,k}\frac{1}{N}\int_{Nt/(t+\varepsilon)}^{N}{\rm d}s\int_{{\mathbb R}^d}{\hat f}({\rm d}\xi)~\prod_{j=1}^d\left(1\wedge\vert s\xi_j\vert^{-2}\right)\notag\\
&\lesssim L_{T,k}\frac{1}{N}\int_{Nt/(t+\varepsilon)}^{N}{\rm d}s\int_{{\mathbb R}^d}{\hat f}({\rm d}\xi)~\left(1\wedge (s\Vert\xi\Vert)^{-2}\right)\notag\\
&\lesssim L_{T,k}\frac{1}{N}\int_{Nt/(t+\varepsilon)}^{N}{\rm d}s\int_{{\mathbb R}^d}{\hat f}({\rm d}\xi)~({s\Vert\xi\Vert})^{-1}
\lesssim L_{T,k}\frac{1}{N}\log
\left(1+\frac{\varepsilon}{t}\right)\notag\\
&\lesssim L_{T,k}\frac{\varepsilon}{Nt}.
\label{4.26}
\end{align}
This, together with (\ref{4.16}) and (\ref{4.25}) implies Lemma \ref{lem4.8}.
\end{proof}
\begin{proof}[Proof of Lemma \ref{lem4.99}]
We can use the same argument as the proof of Lemma 6.2 in \cite{MR4334682} owing
to $f(\R)<\infty$, which we skip here.
\end{proof}
\begin{proof}[Proof of Lemma \ref{lem4.9}]
In this case, recall (\ref{4.17}) and (\ref{4.18}). By using change of variables and semigroup property of the heat kernel we can write
\begin{align}
T_{{\mathcal A}}^2&\leq \frac{L_{T,k}\kappa_{\beta,d}}{N^{\beta}}
\int_0^t\frac{t^\beta}{s^\beta}{\rm d}s\int_{{\mathbb R}^d}{\rm d}\xi~\Vert \xi\Vert^{\beta-d}
\left\vert{\widehat{{\mathcal A}(s)}(t\xi/(Ns))}\right\vert^2,\label{4.27}\\
T_{{\mathcal B}}^2&\leq {L_{T,k}}\int_{t}^{t+\varepsilon}{\rm d}s\int_{{\mathbb R}^{2d}}{\rm d}x{\rm d}y'
\left(I_{N}\ast {\tilde I_{N}}\right)(x)~
{\boldsymbol p}_{2s(t+\varepsilon-s)/(t+\varepsilon)}\left(y'-\frac{s}{t+\varepsilon}x\right)\Vert y' \Vert^{-\beta}.\label{4.28}
\end{align}
Similar to (\ref{4.19}),
\begin{align*}
{\widehat{{\mathcal A}(s)}(t\xi/(Ns))}&=J_1+J_2,
\end{align*}
where
\begin{align*}
J_1&=\int_{[0,1]^d}
\exp\left(\frac{itx\cdot\xi}{t+\varepsilon}\right){\rm d}x~
\exp\left(-\frac{t(t-s)\Vert \xi \Vert^2}{2sN^2}\right)
\left[\exp\left(-\frac{\varepsilon t\Vert \xi \Vert^2}{2(t+\varepsilon)N^2}\right)-1\right],\\
J_2&=\int_{[0,1]^d}
\left[\exp\left(\frac{itx\cdot\xi}{t+\varepsilon}\right)-\exp(ix\cdot\xi)\right]{\rm d}x~
\exp\left(-\frac{t(t-s)\Vert \xi \Vert^2}{2sN^2}\right).
\end{align*}
Therefore,
\begin{align}
T_{{\mathcal A}}^2
&\lesssim \frac{L_{T,k,\beta}}{N^{\beta}}\int_0^t\frac{t^\beta}{s^\beta}{\rm d}s\int_{{\mathbb R}^d}{\rm d}\xi~\Vert \xi\Vert^{\beta-d}
\left\vert J_1\right\vert^2+ \frac{L_{T,k,\beta}}{N^{\beta}}\int_0^t\frac{t^\beta}{s^\beta}{\rm d}s\int_{{\mathbb R}^d}{\rm d}\xi~\Vert \xi\Vert^{\beta-d}
\left\vert J_2\right\vert^2\notag\\
&:={T'}_{{\mathcal A},1}^2+{T'}_{{\mathcal A},2}^2.\label{4.29}
\end{align}
Apply (\ref{4.22}) again to see that
\begin{align}
\left\vert J_1\right\vert^2
&\lesssim\prod_{j=1}^d\phi\left(\frac{t\xi_j}{t+\varepsilon}\right)
\exp\left(-\frac{t(t-s)\Vert \xi \Vert^2}{sN^2}\right)
\left\vert\exp\left(-\frac{\varepsilon t\Vert \xi \Vert^2}{2(t+\varepsilon)N^2}\right)-1\right\vert^2\notag\\
&\lesssim\frac{(t+\varepsilon)^2}{t^2}\frac{1}{\Vert \xi \Vert^2}
\exp\left(-\frac{t(t-s)\Vert \xi \Vert^2}{sN^2}\right)
\left\vert 1-\exp\left(-\frac{\varepsilon t\Vert \xi \Vert^2}{2(t+\varepsilon)N^2}\right)\right\vert^2
\label{4.31}
\end{align}
and
\begin{align}
\left\vert J_2\right\vert^2
%&\lesssim\prod_{j=1}^d
%\left[\left(\frac{\varepsilon^2}{t^2+\varepsilon^2}\right)\wedge\left(\frac{t^2+\varepsilon^2}{(t\xi_j)^2}\right)\right]
%\exp\left(-\frac{t(t-s)\Vert \xi \Vert^2}{sN^2}\right)\notag\\
&\lesssim
\left[\left(\frac{\varepsilon^2}{t^2+\varepsilon^2}\right)\wedge\left(\frac{t^2+\varepsilon^2}{(t\Vert\xi\Vert)^2}\right)\right]
\exp\left(-\frac{t(t-s)\Vert \xi \Vert^2}{sN^2}\right)\notag\\
&\lesssim
\left(\frac{t^2+\varepsilon^2}{t^2}\right)
\left[\left(\frac{\varepsilon^2t^2}{(t^2+\varepsilon^2)^2}\right)\wedge\left(\frac{1}{\Vert\xi\Vert^2}\right)\right].
\label{4.32}
\end{align}
First, we focus on ${T'}_{{\mathcal A},1}^2$. By making the following changes of variables: $z'=\Vert \xi \Vert$, $z=z'/N$, $r=s/t$, we obtain that
\begin{align}
{T'}_{{\mathcal A},1}^2&\lesssim L_{T,k,\beta}\frac{(t+\varepsilon)^2}{tN^2}\int_1^\infty\frac{1}{r^{2-\beta}}\d r\int_0^\infty \frac{z^{\beta-3}}{e^{(t(r-1)z^2)}}\left\vert1-\exp{\left(-\frac{\varepsilon tz^2}{2(t+\varepsilon)}\right)}\right\vert\d z\notag\\
&\lesssim L_{T,k,\beta}\frac{(t+\varepsilon)^2}{tN^2}\int_1^\infty\frac{1}{r^{2-\beta}}\d r\int_0^\infty \frac{z^{\beta-3}}{e^{(t(r-1)z^2)}}\cdot\frac{\varepsilon tz^2}{2(t+\varepsilon)}\d z\notag\\
&\lesssim L_{T,k,\beta}\frac{\varepsilon}{N^2}\int_1^\infty\frac{1}{r^{2-\beta}}\cdot\frac{\Gamma(\beta/2)}{[t(r-1)]^{\beta/2}}\d r\lesssim L_{T,k,\beta}\frac{\varepsilon}{N^2t^{\beta/2}},\label{4.322}
\end{align}
where the last inequality holds by $\int_1^\infty\frac{1}{r^{2-\beta}}\cdot\frac{1}{(r-1)^{\beta/2}}\d r<\infty$ for all $\beta\in(0,2)$.\\
{\rm\bf Case 1: }$\beta\in(0,1)$.
By Lemma \ref{lem6.2} (2), we have
\begin{align}
{T'}_{{\mathcal A},2}^2
&\lesssim
L_{T,k,\b}\frac{(t^2+\varepsilon^2)}{t^2N^{\beta}}
\int_0^t\frac{t^\beta}{s^\beta}{\rm d}s\int_{{\mathbb R}^d}{\rm d}\xi~\Vert \xi\Vert^{\beta-d}
\left[\left(\frac{\varepsilon^2t^2}{(t^2+\varepsilon^2)^2}\right)
\wedge\left(\frac{1}{\Vert\xi\Vert^2}\right)\right]\notag\\
&\lesssim
L_{T,k,\beta}\frac{\varepsilon}{tN^{\beta}}\int_0^t\frac{t^\b}{s^{\beta}}{\rm d}s
\lesssim
\frac{L_{T,k,\beta}}{{N^{\beta}}}\varepsilon.\label{4.34}
\end{align}
Therefore, this together with (\ref{4.322}) concludes that
\begin{align}
{T}_{{\mathcal A}}^2\lesssim L_{T,k}\frac{\e}{(t\wedge1)N^{\beta}}.\label{4.35}
\end{align}
Now we turn to estimate $T_{{\mathcal B}}^2$. Recall (\ref{4.28}). Define $\varphi(x):=\prod_{j=1}^d(1-\vert x_j \vert)$ for all $x\in{\mathbb R}^d$. We observe that ${\left(I_{N}\ast {\tilde I_{N}}\right)}(x)=N^{-d}\varphi(x/N){\bf1}_{[-N,N]^d}(x)$.
 Making a change of variables, it yields that
\begin{align}
{T}_{{\mathcal B}}^2
&\lesssim \frac{L_{T,k}}{N^{\beta}}
\int_t^{t+\varepsilon}{\rm d}s\int_{[-1,1]^d}\varphi(z)\d z~ \E\left(\left\Vert\frac{1}{N}\sqrt{\frac{2s(t+\e-s)}
{t+\e}}Z+\frac{s}{t+\e}z\right\Vert^{-\b}\right)\notag\\
&\lesssim\frac{L_{T,k}}{N^{\beta}}\int_t^{t+\varepsilon}
(t+\e)^{\b}s^{-\b}{\rm d}s\int_{[-1,1]^d}\varphi(z)\Vert z\Vert^{-\b}\d z
\lesssim L_{T,k,\beta}\frac{(t+\e)^{\b}\e}{tN^\b}.
\label{4.36}
\end{align}
Combine (\ref{4.16}) with (\ref{4.35}) and (\ref{4.36}) to conclude (\ref{4.13}).\\
{{\rm\bf Case 2: }}$\beta=1$.
According to Lemma \ref{lem6.4}, we have
\begin{align}
{T'}_{{\mathcal A},2}^2
&\lesssim L_{T,k,\b}\frac{(t^2+\varepsilon^2)}{t^2N}
\notag\\
&~~\times\int_0^t\frac{t}{s}{\rm d}s\int_{{\mathbb R}^d}{\rm d}\xi~\Vert \xi\Vert^{1-d}
\left[\left(\frac{\varepsilon^2t^2}{(t^2+\varepsilon^2)^2}\right)
\wedge\left(\frac{1}{\Vert\xi\Vert^2}\right)\right]\exp\left(-\frac{t(t-s)\Vert \xi \Vert^2}{sN^2}\right)\notag\\
&\lesssim L_{T,k,\b}\frac{(t^2+\e^2)\log N}{Nt}\log_{+}\left(\frac{1}{t}\right)\notag\\
&~~\times\int_{\R^d}{\rm d}\xi~\Vert \xi\Vert^{1-d}\left[\left(\frac{\varepsilon^2t^2}{(t^2+\varepsilon^2)^2}\right)
\wedge\left(\frac{1}{\Vert\xi\Vert^2}\right)\right]\log_{+}\left(\frac{1}{\Vert\xi\Vert}\right).\notag
\end{align}
Then, we made a change of variables $z=\Vert\xi\Vert$ and use Lemma A.3 in \cite{MR4334682} to obtain that
\begin{align}
{T'}_{{\mathcal A},2}^2\lesssim L_{T,k,\b}\frac{\e\log N}{Nt}\log_{+}\left(\frac{1}{t}\right).
\end{align}
Hence, for any $\delta>0$, there exists $L_{T,k,\delta}>0$ such that
\begin{align}
{T}_{{\mathcal A}}^2\lesssim L_{T,k,\delta}\frac{{\e}}{(t\wedge1)^{1+\delta}}.
\label{4.37}
\end{align}
Then it remains to consider the behavior of ${T}_{{\mathcal B}}^2$. Recall (\ref{4.28}). Since the Fourier transform of $I_{N}\ast {\tilde I_{N}}$ is $2^d\prod_{j=1}^d\phi(N\xi_j)$, by using a change of variables $z=Nx$, we obtain
\begin{align*}
{T}_{{\mathcal B}}^2
\leq \frac{{L_{T,k}\kappa_{1,d}}(t+\e)}{(2\pi)^dN}
\int_t^{t+\e}\frac{\d s}{s}\int_{{{\mathbb R}}^d}\prod_{j=1}^d\phi(\xi_j)\exp\left(-\frac{(t+\e)(t+\e-s)\Vert \xi\Vert^2}{sN^2}\right)\Vert \xi\Vert^{1-d}\d \xi,
\end{align*}
where $\phi(x)$ is defined in (\ref{4.30}). Since
$\int_{\R^{d}}\prod_{j=1}^d\phi(\xi_j)\Vert \xi\Vert^{1-d}\d z<\infty$, we conclude that
\begin{align}
{T}_{{\mathcal B}}^2\lesssim L_{T,k}\frac{\e}{Nt}.\label{4.38}
\end{align}
This together with (\ref{4.37}) proves (\ref{4.14}). \\
{{\rm\bf Case 3: }}$\beta\in(1,2\wedge d)$.
Recall (\ref{4.29}) and (\ref{4.32}).
\begin{align}
{T'}_{{\mathcal A},2}^2&\lesssim L_{T,k,\b}\frac{(t^2+\varepsilon^2)}{t^2N^\b}\notag\\
&~~\times\int_0^t\frac{t^\b}{s^\b}{\rm d}s\int_{{\mathbb R}^d}{\rm d}\xi~\Vert \xi\Vert^{\b-d}
\left[\left(\frac{\varepsilon^2t^2}{(t^2+\varepsilon^2)^2}\right)
\wedge\left(\frac{1}{\Vert\xi\Vert^2}\right)\right]\exp\left(-\frac{t(t-s)\Vert \xi \Vert^2}{sN^2}\right)\notag.
\end{align}
Then, by using changes of variables $r'=(t-s)/s$ and $r=-t\Vert\xi\Vert^2r'/N^2$ in order, we obtain that
\begin{align}
{T'}_{{\mathcal A},2}^2&\lesssim
L_{T,k}\frac{(t^2+\varepsilon^2)t^{2-\b}}{t^2N^{2-\beta}}
\int_{{\mathbb R}^d}{\rm d}\xi~\Vert \xi\Vert^{2-\beta-d}
\left[\left(\frac{\varepsilon^2t^2}{(t^2+\varepsilon^2)^2}\right)
\wedge\left(\frac{1}{\Vert\xi\Vert^2}\right)\right]
\int_0^{\infty}r^{\b-2}e^{-r}\d r\notag\\
&\lesssim L_{T,k,\b}\frac{t^{2-\b}\e}{tN^{2-\b}},\label{4.40}
\end{align}
where the last inequality holds by Lemma \ref{lem6.2} (2). This, together with (\ref{4.322}), implies that
\begin{align}
{T}_{{\mathcal A}}^2
\lesssim
L_{T,k,\b}\frac{\e}{(t\wedge1)N^{2-\b}}.
\label{4.41}
\end{align}
Next, we are going to estimate the term ${T}_{{\mathcal B}}^2$. Recall (\ref{4.28}). by making changes of variables we observe that
\begin{align}
{T}_{{\mathcal B}}^2
&\leq \frac{{L_{T,k}\kappa_{\b,d}}(t+\e)^\b}{(2\pi)^dN^\b}
\int_t^{t+\e}\frac{\d s}{s^\b}
\int_{{{{\mathbb R}}}^d}\prod_{j=1}^d\phi(\xi_j)\exp\left(-\frac{(t+\e)(t+\e-s)\Vert \xi\Vert^2}{sN^2}\right)\Vert \xi\Vert^{\b-d}\d \xi\notag\\
&\lesssim
L_{T,k,\b}\frac{t+\e}{N^\b}
\int_{\R^d}\d \xi~\Vert \xi\Vert^{\b-d}\prod_{j=1}^d\phi(\xi_j)
\int_0^{\e/t}\d r~ (1+r)^{\b-2}\exp\left(-\frac{r(t+\e)}{N^2}\Vert \xi\Vert^2\right)\notag\\
&\lesssim
L_{T,k,\b}\frac{(t+\e)\e}{tN^\b}
\int_{\R^d}\Vert \xi\Vert^{\b-d}\prod_{j=1}^d\phi(\xi_j)\d \xi
\lesssim
L_{T,k,\b}\frac{\e}{tN^\b}.\label{4.42}
\end{align}
Finally, (\ref{4.41}) and (\ref{4.42}) together imply (\ref{4.15}).
\end{proof}
Now, we can prove Propositions \ref{prop4.5}-\ref{prop4.7}.
\begin{proof}[Proof of Proposition \ref{prop4.5}]
We assume that $T\geq1$ without losing generality. Choose and fix an arbitrary number $\lambda\in(0,1)$. On the one hand, from Lemma \ref{lem4.8}, for any $\e\in(0,1)$, $N\geq e$ and $t\in(\e^\lambda,T)$ we have
\begin{align}
\Vert S_{N,t+\e}-S_{N,t}\Vert_k
\lesssim L_{T,k,\alpha}\e^{\alpha-\lambda/2}\sqrt{\frac{1}{N}}.\label{4.43}\tag{4.43}
\end{align}
On the other hand, for any $t\in(0,\e^\lambda)$, Lemma \ref{lem4.1} implies that
\begin{align}
\Vert S_{N,t+\e}-S_{N,t}\Vert_k
\leq\Vert S_{N,t+\e}\Vert_k+\Vert S_{N,t}\Vert_k
\lesssim L_{T,k}\e^{\lambda/2}\sqrt{\frac{1}{N}}.\label{4.44}\tag{4.44}
\end{align}
Choose $\lambda=\alpha$ to match the above exponents of $\e$ and define $\gamma:=\alpha/2$ to finish the proof.
\end{proof}
\begin{proof}[Proof of Proposition \ref{prop4.6}]
Similar to the proof of Proposition \ref{prop4.5}, based on Lemma \ref{lem4.2} and \ref{lem4.99},
for any $\delta>0$,
\begin{align}
\Vert S_{N,t+\e}-S_{N,t}\Vert_k
&\lesssim L_{T,k}\e^{(1-\lambda(1+\delta))/2}\sqrt{\frac{\log N}{N}}
~~~{\rm for}~t\in(\e^\lambda,T)
\end{align}
and for any $\alpha\in(0,1)$,
\begin{align}
\Vert S_{N,t+\e}-S_{N,t}\Vert_k
&\lesssim L_{T,k}\e^{\alpha\lambda/2}\sqrt{\frac{\log N}{N}}
~~~{\rm for}~t\in(0,\e^\lambda).
\end{align}
Then, choose $\lambda=1/(1+\alpha+\delta)$, we prove (\ref{4.8}).
\end{proof}
\begin{proof}[Proof of Proposition \ref{prop4.7}]
Similar to the proof of Proposition \ref{prop4.5}, based on Lemma \ref{lem4.3} and \ref{lem4.9}, the following inequalities holds when $\b\in(0,1)$:
\begin{align}
\Vert S_{N,t+\e}-S_{N,t}\Vert_k
&\lesssim L_{T,k,\b}\e^{(1-\lambda)/2}\sqrt{\frac{1}{N^\b}}
~~~{\rm for}~t\in(\e^\lambda,T)\label{4.45}
\end{align}
and
\begin{align}
\Vert S_{N,t+\e}-S_{N,t}\Vert_k
&\lesssim L_{T,k,\beta}\e^{\lambda/2}\sqrt{\frac{1}{N^\b}}
~~~{\rm for}~t\in(0,\e^\lambda).\label{4.46}
\end{align}
Also, if $\b\in(1,2\wedge d)$, we have
\begin{align}
\Vert S_{N,t+\e}-S_{N,t}\Vert_k
&\lesssim L_{T,k,\b}\e^{(1-\lambda)/2}\sqrt{\frac{1}{N^{2-\b}}}
~~~{\rm for}~t\in(\e^\lambda,T)\label{4.47}
\end{align}
and
\begin{align}
\Vert S_{N,t+\e}-S_{N,t}\Vert_k
&\lesssim L_{T,k,\beta}\e^{\lambda(2-\b)/2}\sqrt{\frac{1}{N^{2-\b}}}
~~~{\rm for}~t\in(0,\e^\lambda).\label{4.48}
\end{align}
Therefore, choose $\lambda=1/2$ and $\lambda=1/(3-\b)$ respectively to prove (\ref{4.9}) and (\ref{4.11}).  Finally, based on the proof of Proposition \ref{prop4.6}, we obtain the result in the case that $\b=1$.
\end{proof}
\section{Proofs of Theorems}
In this section, we will establish the weak  convergence of the finite-dimensional distributions to prove the functional CLTs. But first, we need to show a couple of technical lemmas.

Recall (\ref{1.5}). Because of (\ref{1.4}) and a stochastic Fubini argument, we can express $S_{N,t}$ as an It${\rm\hat{o}}$-Walsh stochastic integral:
\begin{align}
S_{N,t}=\int_{\R^{+}\times\R^{d}}\u_{N,t}(s,y)\eta(\d s ,\d y)=\delta(\u_{N,t}),\label{5.1}
\end{align}
where
\begin{align}
\u_{N,t}(s,y)=
{\bf1}_{(0,t)}(s)~U(s,y)\frac{1}{N^d}
\int_{[0,N]^d}{\boldsymbol p}_{s(t-s)/t}{\left(y-\frac{s}{t}x\right){\rm d}x}.\label{5.2}
\end{align}
Since $S_{N,t}$ is Malliavin differentiable (\citealp[Chapter 1.1.3]{MR2200233}), we can write
\begin{align}
D_{s,y}S_{N,t}=
\u_{N,t}(s,y)+\int_{(s,t)\times\R^d}D_{s,y}\u_{N,t}(r,w)\eta(\d r,\d w).\label{5.3}
\end{align}
\begin{lemma}\label{lem5.1}
{\rm($d\geq2$)}. Suppose $\mathcal R(f)<\infty$.
For every $T>0$, we have
\begin{align}
\sup_{t_1,t_2\in(0,T)}\Var\langle DS_{N,t_1},\u_{N,t_2}\rangle_{\mathcal H}\leq L_T{N^{-3}}~~~{\rm for~all~}N\geq e.\label{5.4}
\end{align}
\end{lemma}
\begin{lemma}\label{lem5.2}
{\rm($d=1$)}. Suppose $f(\R)<\infty$,
For every $T>0$, we have
\begin{align}
\sup_{t_1,t_2\in(0,T)}\Var\langle DS_{N,t_1},\u_{N,t_2}\rangle_{\mathcal H}\leq L_T{\left(\frac{N}{\log N}\right)^{-3}}~~~{\rm for~all~}N\geq e.\label{5.5}
\end{align}
\end{lemma}
\begin{lemma}\label{lem5.3}
{\rm(Riesz kernel)}. Suppose $f(\d x)=\Vert x \Vert^{-\b}$, $\b\in(0, 2\wedge d)$.
for any $T>0$ and $t_1,t_2\in(0,T)$, there exists a real number $L=L_{T,t_1,t_2}>0$, such that\\
$\mathbf{(A)}$ If $\beta\in(0,1)$, then
\begin{align}
\Var\langle DS_{N,t_1},\u_{N,t_2}\rangle_{\mathcal H}\leq L{N^{-3\b}}~~~{\rm for~all~}N\geq e.\label{5.6}
\end{align}
$\mathbf{(B)}$ If $\beta=1$, then
\begin{align}
\Var\langle DS_{N,t_1},\u_{N,t_2}\rangle_{\mathcal H}\leq L{\left(\frac{N}{\log N}\right)^{-3}}~~~{\rm for~all~}N\geq e.\label{5.7}
\end{align}
$\mathbf{(C)}$ If $\beta\in(1,2\wedge d)$, then
\begin{align}
\Var\langle DS_{N,t_1},\u_{N,t_2}\rangle_{\mathcal H}\leq L{N^{-3(2-\b)}}~~~{\rm for~all~}N\geq e.\label{5.8}
\end{align}
\end{lemma}
Before proving the above lemmas, we first decompose the quantity $\langle DS_{N,t_1},\u_{N,t_2}\rangle_{\mathcal H}$ by using (\ref{5.3}) and stochastic Fubini argument,
\begin{align*}
&\langle DS_{N,t_1},\u_{N,t_2}\rangle_{\mathcal H}
=\frac{1}{N^{2d}}\int_0^{t_1\wedge t_2}\d s\int_{\R^{2d}}f(\d z)\d y\int_{[0,N]^{2d}}\d x\d x'~U(s,y)U(s,y+z)\\
&~~~~~\times {\boldsymbol p}_{s[(t_1\wedge t_2)-s]/(t_1\wedge t_2)}\left(y-\frac{s}{t_1\wedge t_2}x\right){\boldsymbol p}_{s[(t_1\vee t_2)-s]/(t_1\vee t_2)}\left(y+z-\frac{s}{t_1\vee t_2}x'\right)\\
&~~+\frac{1}{N^{2d}}\int_{[0,t_1]\times\R^d}\eta(\d r,\d w)\int_0^{r\wedge t_2}\d s\int_{\R^{2d}}f(\d z)\d y\int_{[0,N]^{2d}}\d x\d x'~\\
&~~~~~\times {\boldsymbol p}_{r(t_1-r)/t_1}\left(w-\frac{r}{t_1}x\right){\boldsymbol p}_{s(t_2-s)/ t_2}\left(y+z-\frac{s}{t_2}x'\right)U(s,y+z)D_{s,y}U(r,w)\\
&:={\mathcal X}_{N,t_1,t_2}+{\mathcal Y}_{N,t_1,t_2}.
\end{align*}
As a consequence,
\begin{align}
\Var\langle DS_{N,t_1},\u_{N,t_2}\rangle_{\mathcal H}
\leq2\left(\Var{\mathcal X}_{N,t_1,t_2}+\Var{\mathcal Y}_{N,t_1,t_2}\right)
=\frac{2}{N^{4d}}\left(\Phi^{(1)}_{N,t_1,t_2}+\Phi^{(2)}_{N,t_1,t_2}\right),
\label{5.9}
\end{align}
where
\begin{align}
&\Phi^{(1)}_{N,t_1,t_2}
=\int_{[0,t_1\wedge t_2]^2}\d s_1\d s_2\int_{\R^{4d}}f(\d z_1)f(\d z_2)\d y_1\d y_2\int_{[0,N]^{4d}}\d x_1\d x_1'\d x_2\d x_2'\notag\\
&\times {\boldsymbol p}_{s_1[(t_1\wedge t_2)-s_1]/(t_1\wedge t_2)}\left(y_1-\frac{s_1x_1}{t_1\wedge t_2}\right)
{\boldsymbol p}_{s_2[(t_1\wedge t_2)-s_2]/(t_1\wedge t_2)}\left(y_2-\frac{s_2x_2}{t_1\wedge t_2}\right)\notag\\
&\times {\boldsymbol p}_{s_1[(t_1\vee t_2)-s_1]/(t_1\vee t_2)}\left(y_1+z_1-\frac{s_1x_1'}{t_1\vee t_2}\right)
{\boldsymbol p}_{s_2[(t_1\vee t_2)-s_2]/(t_1\vee t_2)}\left(y_2+z_2-\frac{s_2x_2'}{t_1\vee t_2}\right)\notag\\
&\times\Cov[U(s_1,y_1)U(s_1,y_1+z_1),U(s_2,y_2)
U(s_2,y_2+z_2)]\label{5.10}
\end{align}
and
\begin{align}
&\Phi^{(2)}_{N,t_1,t_2}
=\int_{[0,t_1]\times\R^{2d}}\d w f(\d b)\d r
\int_{[0,r\wedge t_2]^2}\d s_1\d s_2
\int_{\R^{4d}}f(\d z_1)f(\d z_2)\d y_1\d y_2\notag\\
&\times\int_{[0,N]^{4d}}\d x_1\d x_1'\d x_2\d x_2'{\boldsymbol p}_{r(t_1-r)/t_1}\left(w-\frac{r}{t_1}x_1\right)
{\boldsymbol p}_{r(t_1-r)/t_1}\left(w+b-\frac{r}{t_1}x_2\right)\notag
\\
&\times {\boldsymbol p}_{s_1(t_2-s_1)/ t_2}\left(y_1+z_1-\frac{s_1}{t_2}x_1'\right){\boldsymbol p}_{s_2(t_2-s_2)/t_2}\left(y_2+z_2-\frac{s_2}{t_2}x_2'
\right)\notag\\
&\times\E[U(s_1,y_1+z_1)U(s_2,y_2+z_2)
D_{s_1,y_1}U(r,w)D_{s_2,y_2}U(r,w+b)].\label{5.11}
\end{align}
Now we turn to estimate the terms $\Phi^{(1)}_{N,t_1,t_2}$ and $\Phi^{(2)}_{N,t_1,t_2}$. By the Poincar\rm{$\acute{\rm e}$}-type inequality (\ref{2.4}), Lemma \ref{lem2.1}, semigroup property and symmetry we can write
\begin{align}
&\Phi^{(1)}_{N,t_1,t_2}
\lesssim
L_T\int_{[0,t_1]^2}\d s_1\d s_2\int_0^{s_1\wedge s_2}\d r\int_{\R^{5d}}f(\d z_1)f(\d z_2)f(\d b)\d y_1\d y_2\notag\\
&\times\int_{[0,N]^{4d}}\d x_1\d x_1'\d x_2\d x_2'
{\boldsymbol p}_{s_1(t_1-s_1)/t_1}\left(y_1-\frac{s_1}{t_1}x_1\right)
{\boldsymbol p}_{s_2(t_1-s_2)/t_1}\left(y_2-\frac{s_2}{t_1}x_2\right)\notag\\
&\times{\boldsymbol p}_{s_1(t_2-s_1)/t_2}\left(y_1+z_1-\frac{s_1}
{t_2}x_1'\right)
{\boldsymbol p}_{s_2(t_2-s_2)/t_2}\left(y_2+z_2-\frac{s_2}{t_2}x_2'\right)\notag\\
&\times{\boldsymbol p}_{r(s_1-r)/s_1+r(s_2-r)/s_2}
\left(b-\frac{r}{s_2}y_2+\frac{r}{s_1}y_1\right),\label{5.12}
\end{align}
where we assume $t_1<t_2$ to simplify the expressions. Similarly,
\begin{align}
&\Phi^{(2)}_{N,t_1,t_2}
\lesssim L_T\int_0^{t_1}\d r\int_{[0,r\wedge t_2]^2}\d s_1\d s_2\int_{\R^{6d}}f(\d b)f(\d z_1)f(\d z_2)\d y_1\d y_2\d w\notag\\
&~~\times\int_{[0,N]^{4d}}\d x_1\d x_1'\d x_2\d x_2'
{\boldsymbol p}_{r(t_1-r)/t_1}\left(w-\frac{r}{t_1}x_1\right)
{\boldsymbol p}_{r(t_1-r)/t_1}\left(w+b-\frac{r}{t_1}x_2\right)\notag\\
&~~\times
{\boldsymbol p}_{s_1(t_2-s_1)/ t_2}\left(y_1+z_1-\frac{s_1}{t_2}x_1'\right)
{\boldsymbol p}_{s_2(t_2-s_2)/t_2}
\left(y_2+z_2-\frac{s_2}{t_2}x_2'\right)\notag\\
&~~\times
{\boldsymbol p}_{s_1(r-s_1)/r}\left(y_1-\frac{s_1}{r}w\right)
{\boldsymbol p}_{s_2(r-s_2)/r}\left(y_2-\frac{s_2}{r}(w+b)\right)\notag\\
&=L_T\int_0^{t_1}\d r\int_{[0,r\wedge t_2]^2}\d s_1\d s_2\int_{\R^{4d}}f(\d b)f(\d z_1)f(\d z_2)\d w
\int_{[0,N]^{4d}}\d x_1\d x_1'\d x_2\d x_2'\notag\\
&~~\times
{\boldsymbol p}_{r(t_1-r)/t_1}\left(w+b-\frac{r}{t_1}x_2\right)
{\boldsymbol p}_{{s_2(t_2-s_2)/t_2}+{s_2(r-s_2)/r}}
\left(z_2-\frac{s_2}{t_2}x_2'+\frac{s_2}{r}(w+b)\right)
\notag\\
&~~\times
{\boldsymbol p}_{r(t_1-r)/t_1}\left(w-\frac{r}{t_1}x_1\right)
{\boldsymbol p}_{{s_1(t_2-s_1)/t_2}+{s_1(r-s_1)/r}}
\left(z_1-\frac{s_1}{t_2}x_1'+\frac{s_1}{r}w\right),
\label{5.13}
\end{align}
where in the first inequality we use the {\rm H${\rm \ddot{o}}$lder} inequality and Lemma \ref{lem2.1}, and in the last equality we use the semigroup property of the heat kernel.

Now we begin to prove Lemmas \ref{lem5.1}-\ref{lem5.3}. Actually, the proof is similar to that of Theorems 1.3-1.5 in \cite{MR4242879}, but there are still slight differences in detail.
\begin{proof}[Proof of Lemma 5.1.]
In order to prove (\ref{5.4}), it suffices to show that
\begin{align}
\sup_{t_1,t_2\in(0,T)}N^{-4d+3}(\Phi^{(1)}_{N,t_1,t_2}
+\Phi^{(2)}_{N,t_1,t_2})\leq L_T~~~{\rm for~all~}N\geq e.\label{5.14}
\end{align}
We first estimate the term $\Phi^{(1)}_{N,t_1,t_2}$, because the time variables $t_1, t_2$ may take different values, we can't use the elementary relation ${\boldsymbol p}_{t}(x){\boldsymbol p}_{t}(y)=2^d\p_{2t}(x+y){\boldsymbol p}_{2t}(x-y)$
as the proof of Theorem 1.3 in \cite{MR4242879}. Therefore, we proceed in the following order: using the changes of variables at first [(i)$m_1=y_1-s_1x_1/t_1$, $m_2=y_2-s_2x_2/t_1 $. (ii)$x_i\mapsto Nx_i$, $x_i'\mapsto Nx_i'$ for $i=1,2$.  (iii)$r_1={s_1N}/{\tau}$, $r_2={s_2N}/{\tau}$, $\sigma=rN/\tau.$], then applying Parseval's identity, to obtain that
\begin{align}
N^{-4d+3}\Phi^{(1)}&_{N,t_1,t_2}
\lesssim L_T\frac{\tau^3}{(2\pi)^{3d}}
\int_{[0,\infty)^3}\d r_1\d r_2\d \sigma
\int_{\R ^3}{\hat f}(\d \xi_1){\hat f}(\d \xi_2){\hat f}(\d \xi_3)\vert \Delta_1(\xi_1,\xi_2,\xi_3)\vert\notag\\
&\times\exp\left(-\frac{\gamma_{1,N}}{2}
\left\Vert\frac{\sigma}{r_1}\xi_3+\xi_1\right\Vert^2
-\frac{\gamma_{2,N}}{2}
\left\Vert\frac{\sigma}{r_2}\xi_3-\xi_2\right\Vert^2
-\frac{\gamma_{3,N}}{2}
\left\Vert\xi_1\right\Vert^2\right.\notag\\
&~~~~~~~~~~~~~~~~~~~~~~~~~~~~~~~~~~~~~~~~~~~~~~-\frac{\gamma_{4,N}}{2}
\left\Vert\xi_2\right\Vert^2
-\frac{\gamma_{5,N}}{2}
\left\Vert\xi_3\right\Vert^2\Bigg)\notag\\
&\lesssim
L_T{\tau^3}
\int_{[0,\infty)^3}\d r_1\d r_2\d \sigma
\int_{\R ^3}{\hat f}(\d \xi_1){\hat f}(\d \xi_2){\hat f}(\d \xi_3)
\vert\Delta_1(\xi_1,\xi_2,\xi_3)\vert,\label{5.15}
\end{align}
where
\begin{align*}
&\Delta_1(\xi_1,\xi_2,\xi_3)
:=\int_{[0,1]^{4d}}\d x_1\d x_1'\d x_2\d x_2'\\
&~~\times\exp[-i(r_1\xi_1+\sigma\xi_3)\cdot (\tau_1x_1)-i(r_2\xi_2-\sigma\xi_3)\cdot (\tau_1x_2)+ir_1\xi_1\cdot(\tau_2x_1')+ir_2\xi_2
\cdot(\tau_2x_2')]
\end{align*}
and
\begin{align*}
\gamma_{1,N}&=\frac{r_1\tau}{N}\left(1-\frac{r_1\tau}{Nt_1}\right),
~~~~~
\gamma_{2,N}=\frac{r_2\tau}{N}\left(1-\frac{r_2\tau}{Nt_1}\right),
~~~~~
\gamma_{3,N}=\frac{r_1\tau}{N}\left(1-\frac{r_1\tau}{Nt_2}\right),
~~~~~\\
\gamma_{4,N}&=\frac{r_2\tau}{N}\left(1-\frac{r_2\tau}{Nt_2}\right),
~~~~~
\gamma_{5,N}=\frac{\sigma\tau}{N}\left(2-\frac{\sigma}{r_1}-\frac{\sigma}{r_2}\right),
~~~~~
\end{align*}
where $\tau_1,\tau_2~{\rm and}~\tau$ are defined in ({\ref{3.1}}). Therefore, we can use the same computations as Theorem 1.3 in \cite{MR4242879} to find that
\begin{align}
N^{-4d+3}&\Phi^{(1)}_{N,t_1,t_2}
\lesssim
L_T{\tau^3}
\int_{[0,\infty)^3}\d r_1\d r_2\d \sigma
\int_{\R ^3}{\hat f}(\d \xi_3){\hat f}(\d \xi_4){\hat f}(\d \xi_5)~
(\Vert\xi_3\Vert\Vert\xi_4\Vert\Vert\xi_5\Vert)^{-1}\notag\\
\times&[1\wedge(\tau_2x)^{-1}][1\wedge(\tau_2y)^{-1}]
[1\wedge(\tau_1\Vert xe_3+ze_5\Vert)^{-1}][1\wedge(\tau_1\Vert ye_4-ze_5\Vert)^{-1}]\notag\\
&\lesssim
L_T\frac{\tau^3}{(\tau_1\wedge\tau_2)^3}
\int_{[0,\infty)^3}\d r_1\d r_2\d \sigma
\int_{\R ^3}{\hat f}(\d \xi_3){\hat f}(\d \xi_4){\hat f}(\d \xi_5)~
(\Vert\xi_3\Vert\Vert\xi_4\Vert\Vert\xi_5\Vert)^{-1}\notag\\
\times&
(1\wedge x^{-1})(1\wedge y^{-1})
[1\wedge(\Vert xe_3+ze_5\Vert)^{-1}][1\wedge(\Vert ye_4-ze_5\Vert)^{-1}]\notag\\
&\lesssim L_T\frac{\tau^3}{(\tau_1\wedge\tau_2)^3}.\label{5.16}
\end{align}
As for $\Phi^{(2)}_{N,t_1,t_2}$, similar to the estimation of $\Phi^{(1)}_{N,t_1,t_2}$, we proceed in the following order:
using the changes of variables at first [(i)$z=w-rx_1/t_1$. (ii)$x_i\mapsto Nx_i,x_i'\mapsto Nx_i'$ for $i=1,2$.  (iii)$r_1={s_1N}/{\tau}$, $r_2={s_2N}/{\tau},\sigma=rN/\tau.$], then applying Parseval's identity and some computations, to obtain that
\begin{align}
N^{-4d+3}\Phi^{(2)}_{N,t_1,t_2}
&\lesssim
L_T{\tau^3}
\int_{[0,\infty)^3}\d r_1\d r_2\d \sigma
\int_{\R ^3}{\hat f}(\d \xi_1){\hat f}(\d \xi_2){\hat f}(\d \xi_3)
\vert\Delta_2(\xi_1,\xi_2,\xi_3)\vert\notag\\
&\lesssim \frac{L_T\tau^3}{(\tau_1\wedge\tau_2)^3},\label{5.17}
\end{align}
where
\begin{align*}
&\Delta_2(\xi_1,\xi_2,\xi_3)
:=\int_{[0,1]^{4d}}\d x_1\d x_1'\d x_2\d x_2'\\
&~~\times\exp[-i(r_1\xi_1+\sigma\xi_2)\cdot (\tau_1x_1)+i(\sigma\xi_2-r_2\xi_3)\cdot (\tau_1x_2)+ir_1\xi_1\cdot(\tau_2x_1')+ir_2\xi_3\cdot(\tau_2x_2')].
\end{align*}
Combining (\ref{5.16}) and (\ref{5.17}), we obtain the desired result.
\end{proof}
\begin{proof}[Proof of Lemma 5.2.]
Recall the decomposition (\ref{5.9}) of $\Var\langle DS_{N,t_1},\u_{N,t_2}\rangle_{\mathcal H}$. In this case, we only need to show that
\begin{align}
\sup_{t_1,t_2\in(0,T)}(\Phi^{(1)}_{N,t_1,t_2}
+\Phi^{(2)}_{N,t_1,t_2})\leq L_TN(\log N)^3~~~{\rm for~all~}N\geq e.\label{5.18}
\end{align}
Recall (\ref{5.12}). One can follow the arguments in the proof of Theorem 1.4 in \cite{MR4242879} to estimate $\Phi^{(1)}_{N,t_1,t_2}$ such that
\begin{align}
\Phi^{(1)}_{N,t_1,t_2}
\lesssim
L_Tt_1t_2^2N\left[2\log N+\log\left(\frac{1}{t_1}+1\right)\right]^3.\label{5.19}
\end{align}
Next, similarly, we discuss $\Phi^{(2)}_{N,t_1,t_2}$ as follows: integrating the variables $x_1$ and $x_2$ on $\R$, applying semigroup property and using the changes of variables and Plancherel's identity, to obtain that
\begin{align}
\Phi^{(2)}_{N,t_1,t_2}
&\lesssim
L_T{Nf(\R)^3}\int_{\R}\phi(z)\d z\int_0^{t_1}\d r\frac{t_1^2t_2}{r}\left(\int_0^{r\wedge t_2}\d s\frac{1}{s}
e^{-\frac{[t_2(t_2-s)/s+t_2^2(r-s)/rs]z^2}{2N^2}}\right)^2\notag\\
&=L_TNf(\R)^3
\int_{\R}\phi(z)\d z\int_0^{t_1\wedge t_2}\d r\frac{t_1^2t_2}{r}\left(\int_0^{r}\d s\frac{1}{s}
e^{-\frac{[t_2(t_2-s)/s+t_2^2(r-s)/rs]z^2}{2N^2}}\right)^2\notag\\
&~~~+L_TNf(\R)^3
\int_{\R}\phi(z)\d z\int_{t_1\wedge t_2}^{t_1}\d r\frac{t_1^2t_2}{r}\left(\int_0^{t_2}\d s\frac{1}{s}
e^{-\frac{[t_2(t_2-s)/s+t_2^2(r-s)/rs]z^2}{2N^2}}\right)^2\notag\\
&:=\Phi^{(2,1)}_{N,t_1,t_2}+\Phi^{(2,2)}_{N,t_1,t_2}.\label{5.20}
\end{align}
With the changes of variables $\theta=(r-s)/s$ and $\xi=(t_2-r)/r$, we have
\begin{align}
\Phi^{(2,1)}_{N,t_1,t_2}
\lesssim L_TNf(\R)^3t_1^2t_2\int_{\R}\phi(z)\d z
\left(\int_0^\infty\frac{1}{\theta+\frac{t_2z^2}{N^2}}e^{-\theta}\d \theta\right)^3
\lesssim L_TN(\log N)^3.\label{5.21}
\end{align}
As for $\Phi^{(2,2)}_{N,t_1,t_2}$, by using the changes of variables $\theta'=(t_2-s)/s$ and $\theta=t_2z^2\theta'/2N^2$ in order, we can write
\begin{align}
\Phi^{(2,2)}_{N,t_1,t_2}
&\lesssim L_TNf(\R)^3\int_{\R}\phi(z)\d z\int_{t_1\wedge t_2}^{t_1}\d r\frac{t_1^2t_2}{r}
\left(\int_0^{\infty}\frac{1}{\theta+1}e^{-\frac{t_2\theta z^2}{2N^2}}\d \theta\right)^2\notag\\
&=L_TNt_1^2t_2\log\left(\frac{t_1}{t_1\wedge t_2}\right)\int_{\R}\phi(z)\d z \left(\int_0^{\infty}\frac{1}{\theta+\frac{t_2z^2}{2N^2}}e^{-\theta}\d \theta\right)^2\notag\\
&\lesssim
L_TN(\log N)^2.\label{5.22}
\end{align}
The proof of Lemma \ref{lem5.2} is now completed.
\end{proof}
\begin{remark}{\rm Comparing with the case $t_1=t_2=t$, we need to deal with the term $\Phi^{(2,2)}_{N,t_1,t_2}$ additionally because of the possibility of the case that $t_2<r<t_1$.}
\end{remark}
\begin{proof}[Proof of Lemma 5.3.]
According to the
proof of (6.14) and (6.21) in \cite{MR4242879}, we first write
\begin{align}
\Phi^{(1)}_{N,t_1,t_2}
&\lesssim
L_TN^{4d-3\b}\int_{[0,t_1]^2}\d s_1\d s_2\int_0^{s_1\wedge s_2}\d r
\left(\frac{\tau}{r}\right)^\b\left(\frac{\tau}{s_1}\right)^\b
\left(\frac{\tau}{s_2}\right)^\b\int_{[0,1]^{4d}}
\d x_1\d x_1'\d x_2\d x_2'\notag\\
&~~~\times \E\left[
\left\Vert\frac{\tau}{Ns_1}\sqrt{\frac{s_1(t_2-s_1)}{t_2}}Z_3
-\frac{\tau}{Ns_2}\sqrt{\frac{s_1(t_1-s_1)}{t_1}}Z_1-
(\tau_1 x_1-\tau_2 x_1')\right\Vert^{-\b}\right.\notag\\
&~~~\times\left\Vert\frac{\tau}{Ns_2}\sqrt{\frac{s_2(t_2-s_2)}{t_2}}Z_4
-\frac{\tau}{Ns_2}\sqrt{\frac{s_2(t_1-s_2)}{t_1}}Z_2-
(\tau_1 x_2-\tau_2 x_2')\right\Vert^{-\b}\notag\\
&~~~\times\left\Vert\frac{\tau}{Nr}\sqrt{\frac{r(s_1-r)}{s_1}+
\frac{r(s_2-r)}{s_2}}Z_5
+\frac{\tau}{Ns_2}\sqrt{\frac{s_2(t_1-s_2)}{t_1}}Z_2\right.\notag\\
&~~~~~~~~~~~~~~~~~~~~~~~~~~~
-\left.\left.\frac{\tau}{Ns_1}\sqrt{\frac{s_1(t_1-s_1)}{t_1}}Z_1
+(\tau_1 x_2-\tau_1 x_1)\right\Vert^{-\b}
\right]\label{5.23}
\end{align}
and
\begin{align}
&\Phi^{(2)}_{N,t_1,t_2}
\lesssim
L_TN^{4d-3\b}\int_0^{t_1}\d r\int_{[0,r\wedge t_2]^2}
\d s_1\d s_2
\left(\frac{\tau}{r}\right)^\b\left(\frac{\tau}{s_1}\right)^\b
\left(\frac{\tau}{s_2}\right)^\b\int_{[0,1]^{4d}}
\d x_1\d x_1'\d x_2\d x_2'\notag\\
&\times \E\left[
\left\Vert\frac{\tau}{Nr}\sqrt{\frac{r(t_1-r)}{t_1}}Z_2
-\frac{\tau}{Nr}\sqrt{\frac{r(t_1-r)}{t_1}}Z_1+
(\tau_1 x_2-\tau_1 x_1)\right\Vert^{-\b}\right.\notag\\
&\times\left\Vert\frac{\tau}{Ns_1}
\sqrt{\frac{s_1(t_2-s_1)}{t_2}+\frac{s_1(r-s_1)}{r}}Z_3
-\frac{\tau}{Nr}\sqrt{\frac{r(t_1-r)}{t_1}}Z_1+
(\tau_2 x_1'-\tau_1 x_1)\right\Vert^{-\b}\notag\\
&\left.\times\left\Vert\frac{\tau}{Ns_2}
\sqrt{\frac{s_2(t_2-s_2)}{t_2}
+\frac{s_2(r-s_2)}{r}}Z_4
-\frac{\tau}{Nr}\sqrt{\frac{r(t_1-r)}{t_1}}Z_2
+(\tau_2 x_2'-\tau_1 x_2)\right\Vert^{-\b}
\right],\label{5.24}
\end{align}
where $Z_1,Z_2,Z_3,Z_4,Z_5$ denote i.i.d.,
$d$-dimension standard normal random variables. Then,
in a similar way to the proof of Theorem 1.5 in \cite{MR4242879},
we conclude that\\
{{\bf Case 1}: $\b\in(0,1)$.}
\begin{align}
\Phi^{(1)}_{N,t_1,t_2}
&\lesssim
L_T\frac{N^{4d-3\b}}{(\tau_1\tau_2^2)^d}
\int_{[0,\tau]^2}\d s_1\d s_2\int_0^{s_1\wedge s_2}\d r
\left(\frac{\tau}{r}\right)^\b\left(\frac{\tau}{s_1}\right)^\b
\left(\frac{\tau}{s_2}\right)^\b
\lesssim L_{T,t_1,t_2}N^{4d-3\b},\label{5.25}\\
\Phi^{(2)}_{N,t_1,t_2}
&\lesssim
L_T\frac{N^{4d-3\b}}{(\tau_1\tau_2^2)^d}
\int_0^{t_1}\d r\int_{[0,r]^2}\d s_1\d s_2
\left(\frac{\tau}{r}\right)^\b\left(\frac{\tau}{s_1}\right)^\b
\left(\frac{\tau}{s_2}\right)^\b
\lesssim L_{T,t_1,t_2}N^{4d-3\b},\label{5.26}
\end{align}
which yield (\ref{5.6}).\\
{{\bf Case 2}: $\b\in(1,2\wedge d)$.}
\begin{align}
&\Phi^{(1)}_{N,t_1,t_2}
\lesssim
L_T\frac{N^{4d-3\b}}{(\tau_1\tau_2^2)^d}
\int_{[0,t_2]^2}\d s_1\d s_2\int_0^{s_1\wedge s_2}\d r\notag\\
&~~~~~\times\left[\left[N^\b\left(\frac{s_2(t_2-s_2)}{t_2}\right)^{-\b/2}
\right]\wedge\left(\frac{t_2}{s_2}\right)^{\b}\right]
\left[\left[N^\b\left(\frac{s_1(t_2-s_1)}{t_2}\right)^{-\b/2}
\right]\wedge\left(\frac{t_2}{s_1}\right)^{\b}\right]\notag\\
&~~~~~\times\left[\left[N^\b\left(\frac{r(s_1-r)}{s_1}+\frac{r(s_2-r)}{s_2}
\right)^{-\b/2}
\right]\wedge\left(\frac{t_2}{r}\right)^{\b}\right]\notag\\
&~~~~\lesssim
L_{T,t_1,t_2}N^{4d-3\b+6},\label{5.27}\\
&\Phi^{(2)}_{N,t_1,t_2}
\lesssim
L_T\frac{N^{4d-3\b}}{(\tau_1\tau_2^2)^d}
\int_0^{t_1}\d r \int_{[0,r]^2}\d s_1\d s_2
\notag\\
&~~~~~\times\left[\left[N^\b\left(\frac{r(t_1-r)}{t_1}\right)^{-\b/2}
\right]\wedge\left(\frac{\tau}{s_1}\right)^{\b}\right]
\left[\left[N^\b\left(\frac{s_1(r-s_1)}{r}\right)^{-\b/2}
\right]\wedge\left(\frac{\tau}{s_1}\right)^{\b}\right]\notag\\
&~~~~~\times\left[\left[N^\b\left(\frac{s_2(r-s_2)}{r}
\right)^{-\b/2}
\right]\wedge\left(\frac{\tau}{s_2}\right)^{\b}\right]\notag\\
&~~~~\lesssim L_T\frac{N^{4d-3\b}}{(\tau_1\tau_2^2)^d}
\int_0^{t_1}\d r \int_{[0,r]^2}\d s_1\d s_2
\left(1\vee\tau^{\b}_1\right)^3
\left[\left[N^\b\left(\frac{r(t_1-r)}{t_1}\right)^{-\b/2}
\right]\wedge\left(\frac{t_1}{s_1}\right)^{\b}\right]\notag\\
&~~~~~\times
\left[\left[N^\b\left(\frac{s_1(r-s_1)}{r}\right)^{-\b/2}
\right]\wedge\left(\frac{t_1}{s_1}\right)^{\b}\right]
\left[\left[N^\b\left(\frac{s_2(r-s_2)}{r}
\right)^{-\b/2}
\right]\wedge\left(\frac{t_1}{s_2}\right)^{\b}\right]\notag\\
&~~~~\lesssim
L_{T,t_1,t_2}N^{4d-3\b+6},\label{5.28}
\end{align}
which yield (\ref{5.8}).\\
{{\bf Case 3}: $\b=1$.} Notice that (\ref{5.27}) and (\ref{5.28})
still hold for $\b=1$, it is easy to see that
\begin{align}
\Phi^{(1)}_{N,t_1,t_2}
\lesssim L_{T,t_1,t_2}N^{4d-3}(\log N)^3~~~{\rm and}~
\Phi^{(2)}_{N,t_1,t_2}
\lesssim L_{T,t_1,t_2}N^{4d-3}(\log N)^3.\label{5.29}
\end{align}
This proves (\ref{5.7}).
\end{proof}

Now, we proceed with the proofs of Theorems \ref{thm1.1}, \ref{thm1.2} and \ref{thm1.3}.
\begin{proof}[Proof of Theorems.]
We first prove Theorem \ref{thm1.1}, and the proofs of Theorem \ref{thm1.2} and \ref{thm1.3}
follow by the same arguments. Choose and fix some $T>0$, by
Lemma \ref{lem4.1} and Proposition \ref{prop4.5}, a standard
application of Kolmogorov's continuity theorem and
the ${\rm Arzel\grave{a}}$-Ascoli theorem ensures that
$\{\sqrt{N}S_{N,\bullet}\}_{N\geq e}$ is a tight net of processes
on $C[0,T]$, so it remains to prove that the finite-dimensional
distributions of the process $t\mapsto\sqrt{N}S_{N,t}$ converge
to those of $\{{\mathcal G}_t\}_{t\in[0,T]}$.

Let us choose and fix some $T>0$ and $m\geq1$ points $t_1,\ldots,t_m\in(0,T)$.
Consider the random vector $F_N=(F_N^{(1)},\ldots,F_N^{(m)})$
 defined by
$$F_N^{(i)}:=\sqrt{N}S_{N,t_i},~~~{\rm for}~i=1,\ldots,m,$$
and let $G=({\mathcal G}_{t_1},\ldots,{\mathcal G}_{t_m})$ be a
centered Gaussian random vector with covariance matrix
 $(g_{t_i,t_j})_{1\leq i,j\leq m}.$ Recall from (\ref{5.1})
to see that $F_N^{(i)}=\delta{(\sqrt{N}\upsilon_{N,t_i})}$, for all
 $i=1,\ldots,m$. Let $V_{N}^{(i)}=\sqrt{N}\upsilon_{N,t_i}$ and
$V_{N}=(V_{N}^{(1)},\ldots,V_{N}^{(m)})$. Lemma \ref{lem2.2} ensures that
\begin{align}
\vert\E(h(F_N))-\E(h(G))\vert\leq\frac{1}{2}\Vert h''\Vert_{\infty}
\sqrt{\sum_{i,j=1}^m\E
\left(\left\vert g_{t_i,t_j}-
\langle DF_{N}^{(i)},V_N^{(j)}\rangle_{\mathcal H}
\right\vert^2\right)},\label{5.30}
\end{align}
for all $h\in C_b^2(\R^m)$. Therefore, it suffices to show that for
any $i,j=1,\ldots,m$,
\begin{align}
\lim_{N\rightarrow \infty}\E\left(\left\vert g_{t_i,t_j}-
\langle DF_{N}^{(i)},V_N^{(j)}\rangle_{\mathcal H}
\right\vert^2\right)=0.\label{5.31}
\end{align}
On one hand, applying Lemma 5.2 in \cite{MR4334682} and
Proposition \ref{prop3.1}, we find that, as $N\rightarrow\infty$,
\begin{align}
\E\langle DF_{N}^{(i)},V_N^{(j)}\rangle_{\mathcal H}=
\Cov(\sqrt{N}S_{N,t_i},\sqrt{N}S_{N,t_j})\rightarrow
g_{t_i,t_j}.\label{5.32}
\end{align}
On the other hand, we can apply Lemma \ref{lem5.1} to see that
\begin{align}
\lim_{N\rightarrow\infty}
\Var\langle DF_{N}^{(i)},V_N^{(j)}\rangle_{\mathcal H}=0.
\label{5.33}
\end{align}
Combining (\ref{5.32}) and (\ref{5.33}), we finish the proof.
\end{proof}


\begin{appendices}

\section{}\label{secA1}

\renewcommand{\thetheorem}{A.\arabic{lemma}}
\newcounter{lemma}
\setcounter{lemma}{1}
\begin{lemma}\label{lem6.1}
(1). Let $\psi_{\tau_1,\tau_2}(z)$ be defined in (\ref{3.13}), then we have
\begin{align*}
\vert \psi_{\tau_1,\tau_2}(z)\vert\lesssim1\wedge\Vert z\Vert^{-2},
\end{align*}
(2). Assume ${\mathcal R}(f)<\infty$, then
\begin{align*}
\int_0^{\infty}\frac{\d s}{s^d}\int_{\R^d}{\hat f}(\d z)
\left\vert\widehat{\left[\left(I_{m}\ast {\tilde I_{n}}\right)
\left(\frac{\bullet}{s}\right)\right]}(z)\right\vert<\infty.
\end{align*}
\end{lemma}

\begin{proof}
(1). By a few computations, it is easy to find that there exists $L_{\tau_2,d}>0$ such that
\begin{align*}
\left\vert\widehat{{I}_{\tau_2}}(z)\right\vert&=\prod_{j=1}^d\sqrt{\frac{2(1-\cos(\tau_2z_j))}{(\tau_2z_j)^2}}
\leq L_{\tau_2,d}\sqrt{\prod_{j=1}^d\left(1\wedge\frac{1}{(z_j)^2}\right)}\\
&\leq L_{\tau_2,d} \sqrt{1\wedge(\max_j z_j)^{-2}}\leq L_{\tau_2,d} \sqrt{1\wedge\Vert z\Vert^{-2}}.
\end{align*}
Similarly,
\begin{align*}
\left\vert\overline{\widehat{{{{ I}}_{\tau_1}}}(z)}\right\vert\leq L_{\tau_1,d}\sqrt{1\wedge\Vert z\Vert^{-2}}.
\end{align*}
Hence,
\begin{align*}
\left\vert \psi_{\tau_1,\tau_2}(z)\right\vert\leq L_{\tau_1,\tau_2,d}\left\vert\widehat{{I}_{\tau_2}}(z)\right\vert \left\vert\overline{\widehat{{{{ I}}_{\tau_1}}}(z)} \right\vert\lesssim1\wedge\Vert z\Vert^{-2}.
\end{align*}
(2). From Lemma 5.9 in \cite{MR4242879}, we know that ${\mathcal R}(f)<\infty$ is equivalent to $\int_{\R^d}\Vert z\Vert^{-1}{\hat f}(\d z)<\infty$. Hence,
\begin{align*}
&\int_0^{\infty}\frac{\d s}{s^d}\int_{\R^d}{\hat f}(\d z)
\left\vert\widehat{\left[\left(I_{m}\ast {\tilde I_{n}}\right)
\left(\frac{\bullet}{s}\right)\right]}(z)\right\vert\\
&\leq(mn)^d(C_{m,n})^d\int_0^{\infty}{\d s}\int_{\R^d}{\hat f}(\d z)
~\prod_{j=1}^{d}\left(1\wedge\frac{1}{(sz_j)^{2}}\right)\\
&\lesssim(C_{m,n})^d
\int_0^{\infty}{\d s}\int_{\R^d}{\hat f}(\d z)
\left(1\wedge\frac{1}{(s\Vert z\Vert)^2}\right)\\
&\lesssim(C_{m,n})^d
\int_0^{\infty}\left(1\wedge\frac{1}{r^2}\right){\d r}
\int_{\R^d}\Vert z\Vert^{-1}{\hat f}(\d z)<\infty.
\end{align*}
\end{proof}

\setcounter{lemma}{2}
\begin{lemma}\label{lem6.2}
(1). For every $\e\in(0,1)$,
$$\int_0^{\infty}\left({s^{-2}}\wedge\e\right)\d s\leq2\sqrt{\e}.$$
(2). For every $\e\in(0,1)$ and $\beta\in(0,1)$, there exists a real number $C_\beta>0$, such that
$$\int_{\R^d}\Vert\xi\Vert^{\beta-d}
\left(\e\wedge\frac{1}{\Vert\xi\Vert^2}\right)\d \xi\leq C_\beta\sqrt{\e}.$$
\end{lemma}

\begin{proof}(1).
\begin{align*}
\int_0^{\infty}\left({s^{-2}}\wedge\e\right)\d s
&\leq\int_0^1\e~\d s+
\e\int_1^\infty\left(\frac{1}{\e s^2}\wedge1\right)\d s\notag\\
&\leq\e+\sqrt{\e}\left(\int_{\sqrt{\e}}^1\d s
+\int_1^\infty\frac{1}{s^2}~\d s\right)=2\sqrt{\e},
\end{align*}
where we use a change of variables $s\mapsto s/\sqrt{\e}$
in the second inequality.\\
(2). Similarly,
\begin{align*}
\int_{\R^d}\Vert\xi\Vert^{\beta-d}
\left(\e\wedge\frac{1}{\Vert\xi\Vert^2}\right)\d \xi
&=\int_0^\infty r^{\beta-1}
\left(\e\wedge\frac{1}{\Vert\xi\Vert^2}\right)\d r\\
&\leq\e\int_0^1 r^{\beta-1}~\d r
+\e\int_1^\infty r^{\beta-1}
\left(1\wedge\frac{1}{(\sqrt{\e}r)^2}\right)\d r\\
&=\frac{\e}{\beta}+(\sqrt{\e})^{2-\beta}
\left(\int_{\sqrt{\e}}^1 r^{\beta-1}~\d r+
\int_1^\infty r^{\beta-3}~\d r
\right)\\
&=\left(\frac{1}{\beta}+\frac{1}{2-\beta}\right)\e^{1-\b/2}
\leq C_\b\sqrt{\e}.
\end{align*}
\end{proof}

\setcounter{lemma}{3}
\begin{lemma}{\rm (\citealp{MR4334682}, Lemma A.1)}\label{lem6.4}
Define
\begin{align*}
G_{N,t}(x):=\frac{t}{\log N}\int_0^t\exp\left(-\frac{(t-s)t}{s}\cdot
\frac{x^2}{N^2}\right)\frac{\d s}{s}~~~~~{ for~all~}N,t>0~and~x\in\R\backslash\{0\}.
\end{align*}
Then, for every $t>0$ and $x\in\R\backslash\{0\}$,
\begin{align*}
\sup_{N\geq e}G_{N,t}(x)\leq
7t\log_{+}(1/t)\log_{+}(1/\vert x \vert),
\end{align*}
where $\log_{+}(w)=\log(e+w)$ for all $w\geq0$. Moreover,
\begin{align*}
\lim_{N\rightarrow\infty}G_{N,t}(x)=2t,~~~~~{for~every}~t>0~
and~x\in\R.
\end{align*}
\end{lemma}



\end{appendices}


\bibliography{sn-bibliography}% common bib file
%%\input new_submit.bbl

\end{document}